\documentclass[11pt]{article}

\makeatletter

\@addtoreset{equation}{section}
\makeatother

\def\be{\begin{equation}}
\def\ee{\end{equation}}

\def\Lq{\ifmmode {\cal L}_q\else ${\cal L}_q$\fi}

\newtheorem{utv}{Proposition}
\newtheorem{cor}{Corollary}

\makeatletter

\@addtoreset{cor}{utv}
\makeatother

\newtheorem{lem}{Lemma}
\newtheorem{rem}{Remark}

\input{amssym.def}
\input{amssym}
\usepackage{eufrak}

\begin{document}

\title{Weyl approach to representation theory of reflection  equation 
algebra}
\author{P.A. Saponov\footnote{saponov@mx.ihep.su}\\
{\small\it Theory Department, Institute for High Energy Physics, 142281
Protvino, Russia}}

\maketitle

\begin{abstract}

The present paper deals with the representation theory of the reflection
equation algebra, connected 
with a Hecke type $R$-matrix. Up to some reasonable additional
conditions the $R$-matrix is arbitrary
(not necessary originated from quantum groups). We suggest a universal
method of constructing finite 
dimensional irreducible non-commutative representations in the framework
of the Weyl approach well known 
in the representation theory of classical Lie groups and algebras. With
this method  a series of 
irreducible modules is constructed which are parametrized  by Young
diagrams. The spectrum
of central elements $s_k=Tr_q\,L^k$ is calculated in the single-row and
single-column representations.
A rule for the decomposition of the tensor product of modules into the
direct sum of irreducible 
components is also suggested.
\end{abstract}

\section{Reflection Equation Algebra}

Reflection equation and the corresponding algebra which will be called
the reflection equation
algebra (REA for short) play a significant role in the theory of
integrable systems and non-commutative
geometry. In application to integrable systems the reflection equation
with a spectral parameter is 
mainly used. First it appears in the work by I. Cherednik \cite{Cher}.
Usually it comprises the information
about the behaviour of a system at a boundary, for example, describes
the reflection of particles on a
boundary of the configuration space.

The reflection equation without a spectral parameter is important for
the non-commutative geometry.
One of the first applications the corresponding REA found in the theory
of differential calculus on
quantum groups (see, e.g., \cite{IP}). In such a differential calculus
REA with the Hecke type $R$-matrix 
is a non-commutative analog of the algebra of vector fields on the
groups $GL(N)$ or $SL(N)$. 
Besides, REA serves as a base for a definition of quantum analogs of
homogeneous spaces --- 
orbits of the coadjoint representation of a Lie group, as well as
quantum analogs of linear bundles
over such orbits (see, e.g., \cite{DM,GS}). 

In this paper we turn to problems of the representation theory of REA
without a spectral parameter.
We are interested in the following main topics:\par
{\bf i)} a construction of finite dimensional non-commutative
irreducible representations and the
calculation of spectrum (characters) of central elements in these
representations;\par
{\bf ii)} a rule for the decomposition of the tensor product of
irreducible modules
into irreducible components.

Before reviewing the known results, we introduce some necessary
definitions and 
notations.

Consider an associative algebra \Lq{} with the unity $e_{\cal L}$ over
the complex field $\Bbb C$
generated by $n^2$ elements $\hat l_i^{\,j}$,  $1\le i,j\le n$, $n$
being a fixed positive integer.
Let the generators satisfy the following quadratic commutation relations
\be
R_{12}\hat L_1R_{12}\hat L_1 - \hat L_1R_{12}\hat L_1R_{12} = 0, 
\quad \hat L_1\equiv \hat L\otimes I,\label{RE}
\ee
where the matrix  $\hat L\in {\rm Mat}_n(\Lq)$ is composed of $\hat
l_i^{\,j}$: $\hat L
=\|\hat l_i^{\,j}\|$. Here the lower index enumerates the rows while the
upper one columns.
In (\ref{RE}) and  everywhere below the use is made of the compact
matrix notations \cite{FRT}
when the index of an object indicates the vector space to which the
object belongs (or in which this 
object acts). The symbol $I$ stands for the unity matrix whose dimension
is always clear from the context
of formulae. A numerical $n^2\times n^2$ matrix $R$ is a solution of
{\it the Yang-Baxter equation}
\be 
R_{12}R_{23}R_{12} = R_{23}R_{12}R_{23}.\label{YBE}
\ee
The algebra \Lq{} described above will be called {\it the reflection
equation algebra} (REA).

Impose now several additional conditions on the matrix $R$. First of
them is {\it the Hecke condition}
\be
(R-qI)(R+q^{-1}I) = 0.\label{Hec}
\ee
The parameter $q$ is a fixed nonzero complex number with the only 
constraint\footnote{All our subsequent constructions possess a well
defined 
``classical limit'' $q\rightarrow 1$. This limit corresponds to REA with
an involutive 
 $R$-matrix: $R^2=I$.}
\be
q^k\not=1, \quad\forall\,k\in{\Bbb N}.\label{generic}
\ee
As a consequence, the $q$-analogs of all integers are nonzero
\be
k_q\equiv\frac{q^k-q^{-k}}{q-q^{-1}}\not=0,\quad \forall\, k\in {\Bbb
N}.
\label{q-num}
\ee

Besides, we shall suppose the $R$-matrix to be {\it skew-invertible}
that is there 
exists an $n^2\times n^2$ matrix $\Psi$ such that 
$$
\sum_{a,b}R_{ia}^{\;jb}\Psi_{bk}^{\;as} = \delta_i^{\,s}\delta_k^{\,j} = 
\sum_{a,b} \Psi_{ia}^{\;jb}R_{bk}^{\;as}.
$$
In the compact notations the above formula reads
\be
Tr_{(2)}R_{12}\Psi_{23} = P_{13} = Tr_{(2)}\Psi_{12}R_{23},
\label{closed}
\ee
where the symbol $Tr_{(2)}$ means the calculation of trace in the second
space and
$P$ is the permutation matrix.

To formulate the last requirement on $R$ one should consider the
connection of the $A_k$ series
Hecke algebras with the group algebras of finite symmetric groups. One
of the simplest definition
of the Hecke algebra reads as follows.

Fix a nonzero complex number $q$. The {\it Hecke algebra of $A_k$
series} ($k\ge 2$) is an 
associative algebra $H_k(q)$ over the complex field $\Bbb C$ generated
by the unit element 
$1_{H}$ and $k-1$ generators $\sigma_i$ subject to the following
relations:
$$
\left.
\begin{array}{l}
\sigma_i\sigma_{i+1}\sigma_i = \sigma_{i+1}\sigma_i\sigma_{i+1}\\
\rule{0pt}{5mm}
\sigma_i\sigma_j = \sigma_j\sigma_i\qquad\qquad {\rm if\ }|i-j|\ge 2\\
\rule{0pt}{5mm}
(\sigma_i-q\,1_{H})(\sigma_i+q^{-1}\,1_{H}) = 0\\
\end{array}
\right\}\;\; i=1,2,\dots k-1.
$$
In some cases it proves to be convenient to consider $q$ as a formal
parameter and consider
the Hecke algebra over the field of rational functions in the
indeterminate $q$. We shall always
bear in mind this extension when considering the classical limit
$q\rightarrow 1$.

Let us treat $R$ as the matrix of a linear operator (in a fixed basis)
which acts
in the tensor square $V^{\otimes 2}$ of a finite dimensional vector
space $V$,\  $\dim V = n$.
Then an arbitrary Hecke $R$-matrix define {\it the local representation}
of $H_k(q)$ in 
$V^{\otimes k}$
\be 
\sigma_i\rightarrow \rho_R(\sigma_i)= R_{ii+1} = 
I^{\otimes(i-1)}\otimes R\otimes I^{\otimes(k-i-1)}\in {\rm
End}(V^{\otimes k}).
\label{loc-rep} 
\ee

If the parameter $q$ satisfies (\ref{generic}), then for any positive
integer $k$ the Hecke algebra 
$H_k(q)$ is known to be isomorphic to the group algebra ${\Bbb C}[{\cal
S}_k]$ of the $k$-th order
permutation group ${\cal S}_k$. As a consequence, there exist elements
${\cal Y}_{\nu(a)}(\sigma)
\in H_k(q)$ which are the $q$-analogs of the Young idempotents
(projectors) widely used in the theory of 
symmetric groups. These $q$-idempotents are parametrized by standard
Young tableaux $\nu(a)$ corresponding
to each diagram or equivalently to each partition $\nu\vdash k$. The
number of all standard tableaux 
$\nu(a)$ which one can construct for a given $\nu$ will be denoted
$\dim[\nu]$
$$
\dim[\nu] = \#\{\nu(a)\}.
$$
In the local representation (\ref{loc-rep}) the elements ${\cal
Y}_\nu(\sigma)$ are realized as some 
projector operators in $V^{\otimes k}$. 
With respect to the action of these projectors the space $V^{\otimes k}$
is decomposed
into the direct sum of subspaces $V_{\nu}$ as in the case of the
symmetric group
\be
V^{\otimes k} = \bigoplus_{\nu\vdash k}\bigoplus_{a=1}^{\dim[\nu]}
V_{\nu(a)},\quad V_{\nu(a)} = Y_{\nu(a)}(R)
\triangleright V^{\otimes k}.
\label{razl-prostr}
\ee
The projector $Y_{\nu(a)}(R) = \rho_R({\cal Y}_{\nu(a)})$ is given by
some polynomial in matrices $R_{ii+1}$.
For detailed treatment of these questions, explicit formulae for
$q$-projectors and the extensive list 
of original papers the reader is referred to \cite{OgP}. 

So, we shall  assume that there exists an integer $p>0$ such that the
image of the 
$q$-an\-ti\-sym\-met\-ri\-zer 
${\cal A}^{(p+1)}(\sigma)\in H_k(q)$ ($\forall\, k>p$) under the local
$R$-matrix representation 
$\rho_R$ is identical zero while the image of the $q$-antisymmetrizer
${\cal A}^{(p)}(\sigma)\in H_k(q)$ 
is a unit rank projector in the space $V^{\otimes k}$ 
\be
\exists\,p\in{\Bbb N}:\quad 
\left\{
\begin{array}{l}
{\cal A}^{(p+1)}(\sigma)\stackrel{\rho_R}{\longrightarrow} 
A^{(p+1)}(R)\equiv 0,\\
\rule{0pt}{6mm}
{\cal A}^{(p)}(\sigma)\stackrel{\rho_R}{\longrightarrow} 
A^{(p)}(R),\quad
{\rm rank}\,A^{(p)}(R) = 1.
\end{array}
\right.
\label{f-rank}
\ee
Such a number $p$ will be called {\it the symmetry rank} of the matrix
$R$. For example, the symmetry
rank of $R$-matrix connected with the quantum universal enveloping
algebra $U_q(sl_n)$ is equal to $n$. 
Examples of $n^2\times n^2$ $R$-matrices with $p<n$ (for $n\geq 3$) were
found in \cite{Gur}.

Introduce now two $n\times n$ matrices $B$ and $C$
\be
B_1 = Tr_{(2)}\Psi_{21}, \quad C_1 = Tr_{(2)}\Psi_{12},
\label{BC}
\ee
where $\Psi$ is defined in (\ref{closed}). If the $R$-matrix has the
symmetry rank $p$ these matrices 
are nonsingular and their product is a multiple of the unit matrix 
\cite{Gur}
\be 
B\cdot C = \frac{1}{q^{2p}}\,I.
\label{nonsing}
\ee
Besides, $B$ and $C$ have the following traces
\be
TrB=TrC=\frac{p_q}{q^p}.
\label{BC-norm}
\ee
The matrices  $B$ and $C$ play the central role in what follows.

The simplest example of REA is obtained by choosing the $U_q(sl_2)$
$R$-matrix ($n=2$)
$$
R = \left(
\matrix{
q&0&0&0\cr
0&\lambda&1&0\cr
0&1&0&0\cr
0&0&0&q}\right)
\quad \lambda\equiv q - q^{-1},\qquad
\hat L=\left(\matrix{\hat a&\hat b\cr \hat c&\hat d}\right).
$$
In this case equation (\ref{RE}) leads to six permutation relations for
the generators of REA
\be
\begin{array}{r@{\hspace{20mm}}l}
q^2 \hat a\hat b = \hat b\hat a  
&q(\hat  b\hat c - \hat c\hat b)  =\lambda \,\hat a(\hat d-\hat a)\\
q^2 \hat c\hat a = \hat a\hat c  &q(\hat c\hat d - \hat d\hat c)  
= \lambda\, \hat c \hat a\\
\hat a\hat d=\hat d\hat a &  q(\hat d\hat b - \hat b\hat d) = 
\lambda\, \hat a \hat b.
\end{array}
\label{ex:rea}
\ee

Consider a map $Tr_q:{\rm Mat}_n(\Lq)\rightarrow \Lq$ which is called
{\it the quantum trace}
\cite{FRT}
\be
Tr_q(X)\stackrel{\mbox{\tiny def}}{=} Tr(C\cdot X),\quad X\in {\rm
Mat}_n(\Lq).
\label{q-sled}
\ee
One can show that the quantities 
\be
s_m(\hat L) = Tr_q(\,\hat L^m)\quad 1\le m\le p-1
\label{centr-el}
\ee
are {\it independent} central elements of REA (see \cite{FRT}).
Presumably these elements (together with 
$e_{\cal L}$) generate the whole center of REA, but we do not know the
proof of this hypothesis.
Calculation of spectrum of central elements  $s_m$ in irreducible
representations  of \Lq{} is one of our 
aims.

At present there exist rather lot of works, devoted to the
representation theory of REA.
First of all, this algebra possesses a large number of one-dimensional
(commutative) representations.
For example, it is evident that at any choice of $R$-matrix  relation
(\ref{RE}) will be satisfied
if one sets $\hat l_i^{\,j} = \alpha\,\delta_i^{\,j}$.
Less trivial representation can be obtained for our simple example
(\ref{ex:rea}) by putting
$\hat a =0$. Then the remaining generators are represented by three
arbitrary complex numbers.
Such like representations were considered in detail in \cite{KulShweib}
for the REA with
$R$-matrices coming from $U_q(sl_n)$ and its super-symmetric
generalizations.

But since REA is a non-commutative algebra, the images of a part of its
generators
are inevitably zero in any one-dimensional representation. That is the
kernel of any
one-dimensional representation of REA must contain some of its
generators. These representations
are not comprised by our approach and we shall not consider them.

The main object of our interest will be the non-commutative
representations which for {\it any}
generator of REA put into correspondence a nontrivial linear operator in
a finite dimensional vector 
space. An example of such a representation can be constructed in the
following way. Let 
us use the fact that REA (\ref{RE}) is an adjoint comodule over some
Hopf algebra which is similar to 
to the algebra of functions over the quantum group. Suppose the
commutations among the
generators $t_i^{\,j}$ of the Hopf algebra to be given by the matrix
relation \cite{FRT}
\be
R_{12}T_1T_2 = T_1T_2R_{12}.\label{RTT}
\ee
The above multiplication is compatible with the comultiplication
$\Delta$
\be
t_i^{\,j}\,\stackrel{\Delta}{\longrightarrow}\, \sum_k t_i^{\,k}\otimes
t_k^{\,j}.
\label{coprod}
\ee
Then (\ref{RE}) is covariant with respect to the transformation
\cite{KulS}
\be
\hat l_i^{\,j}\rightarrow t_i^{\,k}S(t_p^{\,j})\otimes \hat l_k^{\,p},
\label{tsl}
\ee
where $S(\cdot)$ stands for the antipodal map\footnote{Note, that it is
skew-invertibility
(\ref{closed}) which allows one to define the antipode in (bi)algebra
(\ref{RTT}) 
(to be more precise, in some its extension, see  \cite{FRT}).}. So, if
one knows the representations of
Hopf algebra (\ref{RTT}) then, given a representation of REA, one can
construct another one
on the base of (\ref{tsl}).

Moreover, in the case of $U_q(sl_n)$ $R$-matrix one has an additional
possibility of constructing
non-commutative representations. The matter is that in such a case there
exists an embedding
of REA into $U_q(gl_n)$. At the level of generators this embedding is
described by the formula
\cite{FRT}
\be
\hat L = S(L^-)L^+,
\label{mapinq}
\ee
where $L^\pm$ are matrices composed of the $U_q(gl_n)$ generators.
Therefore, starting from a 
representation of the quantum group one can find the corresponding
representation of REA by means
of (\ref{mapinq}). In recent paper \cite{DoKulMu} this approach was
extended to the case of an
arbitrary quasitriangular Hopf algebra. The authors of the cited paper
constructed a universal 
solution of (\ref{RE}) basing on the universal $R$-matrix of a
quasitriangular Hopf algebra 
$\cal H$. The generators $\hat l_i^{\,j}$ turn out to be elements of the
tensor product of 
$\cal H$ and its  ``twisted dual'' algebra.

However, it is worth pointing out that all the methods mentioned above
are essentially
based on the representation theory of objects which are external to REA,
namely, on the theory
of quasitriangular Hopf algebras, the quantum groups being a particular
case of them.
But as was shown in \cite{BelDr}, the Yang-Baxter equation possesses a
lot of solutions
do not connected with a quantum group (see, also, \cite{Gur}). For such
type solutions
we cannot use map (\ref{mapinq}) and cannot construct the REA
representations on the 
base of quantum group ones. 

Moreover, having fixed the $R$-matrix in (\ref{RE}), one completely
defines all properties 
of REA and its representation theory as well. The situation is similar
to the Lie algebra 
theory, where the set of structure constants defines all properties of
the algebra. 
Therefore, it is quite natural to develop the representation theory of
REA using only 
the given $R$-matrix, that is entirely in terms of REA itself.

The most efficient method to solve this problem seems to consist in the
direct analysis
of the explicit commutation relations of the REA generators. For
$U_q(sl_n)$ $R$-matrix
(at small values of $n$) one can proceed in an analogy with the
representation theory of 
the universal enveloping algebra of a (simple) Lie algebra. Following
this way, P.P. Kulish 
\cite{Kul} succeeded in finding all highest vector representations of
the simplest REA 
(\ref{ex:rea}). Besides two one-dimensional representations this algebra
has a series of 
finite dimensional irreducible non-commutative representations and one
infinite dimensional 
representation --- an analog of the Verma module of the universal
enveloping algebra.

Unfortunately, this approach is not universal. It is in essential
dependence on the particular 
choice of $R$-matrix. The explicit components of matrix relations
(\ref{RE}) may become 
completely different when we change the $R$-matrix. Therefore, one would
have 
to repeat the analysis of commutation relations from the very beginning
for each possible 
$R$-matrix. Another obstacle in this way is more technical. The matter
is that even in the case 
of $U_q(sl_n)$ $R$-matrix the complexity of the explicit form of
(\ref{RE}) increases very 
quickly with growing of $n$. This leads to additional difficulties as
compared with the case 
of universal enveloping algebra, when the commutation relations among
generators can be written 
in a compact form for an arbitrary $n$.

In the present paper we suggest a  universal method of constructing
finite dimensional 
representations of REA generated by (\ref{RE}). These representations
are parametrized 
by Young diagrams and exist for {\it any} $R$-matrix satisfying the
additional conditions 
(\ref{Hec}), (\ref{closed}) and (\ref{f-rank}). For the representations
corresponding to 
single-row and single-column diagrams we calculate the spectrum of
central elements 
(\ref{centr-el}). In the particular case of $U_q(sl_2)$ $R$-matrix our
result reproduces 
the series of finite dimensional representations of REA (\ref{ex:rea})
obtained in \cite{Kul}.

The paper is organized as follows. In Section \ref{t} we construct an
irreducible representation 
of REA with an arbitrary Hecke $R$-matrix possessing a finite symmetry
rank. This representation 
is  called the fundamental one (of $B$ type) since its tensor products
are decomposed into 
irreducible components similarly to those of fundamental vector
representation of $U(gl_n)$. 

In Section \ref{th} we study the $k$-th tensor power of the fundamental
module of $B$ type 
and consider its decomposition into higher dimensional REA modules.

Section \ref{f} is devoted to another fundamental module (of $R$ type).
The construction
is based on the general theory of dual Hopf algebras and can be easily
generalized to the case
of an arbitrary Hecke $R$-matrix. The connection of $B$ and $R$ type
fundamental modules is 
established. We also consider an example of reducible indecomposable
module over REA which is 
not equivalent to either $B$ or $R$ type module. At the end of the
section we give a short 
resum\'e of the obtained results and mention some open questions of the
suggested approach.
\vskip 4mm

\noindent{\bf Acknowledgment} The author would like to thank D.I.
Gurevich and 
P.N. Pyatov for valuable and helpful discussions on the problems
considered in 
this paper. Also I express my gratitude to the staff and guests of the 
Max-Planck Institute f\"ur Mathematik for kind hospitality, perfect
conditions 
for work and encouraging and friendly atmosphere.

\section{Fundamental module of $B$ type}
\label{t}

Consider the REA generated by relations (\ref{RE}) and make the linear
shift of generators
\be
l_i^{\,j} =  \hat l_i^{\,j} + \frac{1}{\lambda}\delta_i^{\,j}\,e_{\cal
L},
\label{shift}
\ee
where $e_{\cal L}$ is the unit element of \Lq{} and $\lambda =
q-q^{-1}$. On taking into account 
the Hecke condition (\ref{Hec}) one obtains the commutation relations
for the new generators
\be
R_{12}L_1R_{12}L_1 - L_1R_{12}L_1R_{12} = R_{12}L_1 - L_1R_{12}.
\label{mREA}
\ee
In what follows we shall call the algebra generated by (\ref{mREA})
{\it the modified reflection equation algebra} (mREA) and retain the
notation  
$\Lq$ for it. Note, that unless $q=1$ mREA is isomorphic to REA with
relations
(\ref{RE}). As a consequence, any representation of mREA can be
transformed into that of 
REA and vice versa. Nevertheless, these algebras are different at the
classical limit
since isomorphism (\ref{shift}) is broken at $q\rightarrow 1$ in virtue
of singularity 
of $\lambda^{-1}$.

In the particular case of $U_q(sl_n)$ $R$-matrix the classical limit of
(\ref{mREA})
gives the commutation relations of the $U(gl_n)$ generators.

So, consider the mREA $\Lq$ generated by the unit element and $n^2$
generators $l_i^{\,j}$ 
with commutation relations (\ref{mREA}). Let us take an $n$-dimensional
vector space $V$ and 
fix an arbitrary basis of $n$ vectors $e_i$, $1\le i\le n$. Define a
{\it linear}
map $\pi:\,\Lq\rightarrow {\rm End}(V)$ in accordance with the rules
\be
\begin{array}{l}
\pi(e_{\cal L}) = {\rm id}_V\\
\rule{0pt}{5mm}
\pi(l_i^{\,j})\triangleright e_k = e_iB_k^{\,j},\label{B-rep}\\
\rule{0pt}{5mm}
\pi(l_1\cdot l_2\cdot \dots\cdot l_k) = \pi(l_1)\cdot
\pi(l_2)\cdot \dots \cdot\pi(l_k), \quad \forall\,k\in{\Bbb N},
\end{array}
\label{B-fun}
\ee
where the matrix $B$ is defined in (\ref{BC}) and ${\rm id}_V$ is the
identical 
operator on $V$.

\begin{utv}\label{u1}  {\bf (\cite{GLS2})}
The linear map (\ref{B-rep}) defines an irreducible representation of
\Lq{}
(\ref{mREA}) in the space $V$. This will be called the fundamental
module of $B$ type.
\end{utv}

\noindent{\bf  Proof}\ \ 
To prove that $\pi$ realizes the representation of $\Lq$ one should only
verify that operators 
(\ref{B-rep}) do satisfy (\ref{mREA}). This can be easily done by a
straightforward calculation.
The only fact needed in this way consists in the following simple
consequence of (\ref{closed}) and 
(\ref{BC})
\be
Tr_{(1)}B_1R_{12} = I.
\label{B-R}
\ee

Irreducibility follows from the non-singularity of $B$ (\ref{nonsing}).
Using this fact one can show 
that the operators $\pi(l_i^{\,j})$ span  ${\rm End}(V)$ and, therefore,
the space $V$ does not
contain  proper invariant subspaces with respect to $\pi$.  \hfill
\rule{6.5pt}{6.5pt}
\medskip

If we make shift (\ref{shift}) in  our example (\ref{ex:rea}) we get 
the commutation relations of
the corresponding mREA 
\be
\begin{array}{l@{\hspace{20mm}}l}
q^2 ab - ba =q b  &q( bc - cb)  =(\lambda a -1)(d-a)\\
q^2 ca - ac = qc  &q(cd - dc)  =  c(\lambda a -1)\\
ad=da &  q(db - bd) = (\lambda a -1) b.
\end{array}
\label{ex:mrea}
\ee
Here in the right column the unit stands for $e_{\cal L}$.
The matrices $B$ and $C$ have the form
$$
B = \left(
\matrix{
q^{-1}&0\cr
0&q^{-3}}
\right),\qquad
C = \left(
\matrix{
q^{-3}&0\cr
0&q^{-1}}
\right).
$$
The fundamental representation (\ref{B-rep}) reads 
\be
\pi(a) = \left(
\matrix{
q^{-1}&0\cr
0&0}
\right),\;\;
\pi(b) = \left(
\matrix{
0&q^{-3}\cr
0&0}
\right),\;\;
\pi(c) = \left(
\matrix{
0&0\cr
q^{-1}&0}
\right),\;\;
\pi(d) = \left(
\matrix{
0&0\cr
0&q^{-3}}
\right).
\label{ex:rep}
\ee

Given a representation of \Lq{} one can find the corresponding
representation of 
the quotient algebra 
\be
{\cal SL}_q = \Lq/\{Tr_qL\},
\label{sl-red}
\ee
where $\{X\}$ stands for the ideal generated by a given subset $X\subset
\Lq$.
The commutation relations among the generators $f_i^{\,j}$ of ${\cal
SL}_q$ has 
the same form (\ref{mREA}) as those of $\Lq$ (with substitution
$L\rightarrow F$
where $F=\|f_i^{\,j}\|$), but now the generators are linear dependent
due to $Tr_qF =0$.

At $q\rightarrow 1$ in the case of $U_q(sl_n)$ $R$-matrix the
commutation relations of 
the ${\cal SL}_q$ generators transform into those of $U(sl_n)$
generators.
For this reason the passage from \Lq{} to ${\cal SL}_q$ (or from the
\Lq{} representation 
to the corresponding ${\cal SL}_q$ one) will be loosely called the
sl-re\-duc\-tion in what 
follows.

The transformation of an irreducible \Lq{} representation $\rho$ acting
in a finite dimensional
space $V$ into the ${\cal SL}_q$ representation $\bar\rho$ is realized
as follows. Due to $Tr_qL$ 
is a central element of \Lq{} and $\rho$ is an irreducible
representation one gets 
$$
\rho(Tr_qL) = \chi(Tr_qL)\,{\rm id}_V\equiv \chi_1\,{\rm id}_V,
$$
where $\chi: Z(\Lq)\rightarrow {\Bbb C}$ is a character of the center
$Z(\Lq)$.
Then the straightforward calculation shows that the ${\cal SL}_q$
generators $f_i^{\,j}$
in representation $\bar\rho$ are given by
\be
\bar\rho(f_i^{\,j}) = \frac{1}{\omega}\,\Bigl(\rho(l_i^{\,j}) - 
\delta_i^{\,j}\,\frac{\chi_1}{TrC}\,{\rm id}_V\Bigr), \quad 
\omega = 1- \lambda\,\frac{\chi_1}{TrC}.
\label{sl-red-rep}
\ee
The traceless property $\bar\rho(Tr_qF) = 0$ is evident and the factor
$\omega^{-1}$ 
ensures the correct normalization of the right hand side of
(\ref{mREA}).
\begin{rem}\label{r1}
As can be easily seen from definition (\ref{RE}) the REA admits the
``renormalization''
automorphism $\hat l_i^{\,j}\rightarrow z \hat l_i^{\,j}$ with nonzero
complex number $z$.
The same is true for the REA representations as well. At the level of
mREA representations
this automorphism reads
\be
\rho(l_i^{\,j})\rightarrow \rho_z(l_i^{\,j}) = z\rho(l_i^{\,j}) +
\delta_i^{\,j}\,\frac{1-z}{\lambda}\,{\rm id}_V,
\label{ren-aut}
\ee
where $\rho$ is an arbitrary mREA representation in the space $V$.
Basing on (\ref{sl-red-rep})
one can show that the corresponding ${\cal SL}_q$ representation $\bar
\rho$ does not depend on 
$z$ that is the whole class of mREA representations $\rho_z$ connected
by renormalization 
automorphism (\ref{ren-aut}) gives the same ${\cal SL}_q$ representation
$\bar \rho$.
\end{rem}
Let us now obtain the sl-reduction of the $B$ type representation $\pi$
defined by (\ref{B-rep}).
In virtue of (\ref{nonsing}) one finds
$$
\chi_1 = \chi(Tr_qL) =  q^{-2p}.
$$
Then, taking into account (\ref{BC-norm}) and (\ref{sl-red-rep}) we find
the $B$ type representation
$\bar \pi$ of algebra (\ref{sl-red})
\be
\bar\pi(f_i^{\,j}) = \frac{1}{\omega}\,\Bigl(\pi(l_i^{\,j}) - 
\frac{\delta_i^{\,j}}{q^pp_q}\,{\rm id}_V\Bigr), \quad 
\omega = \frac{q^{1-p}}{p_q}\,(q^{p-2}(p+1)_q - 1).
\label{rep-red}
\ee

Consider again our example (\ref{ex:mrea}) of mREA \Lq. To get the
corresponding algebra 
${\cal SL}_q$ it is necessary to take the quotient of $\Lq$ over the
ideal generated by
$Tr_qL$. Using the explicit form of $C$ one has
$$
Tr_qL = \frac{1}{q^3}\,a+\frac{1}{q}\,d.
$$
With a new generator $h=a-d$ one can rewrite the commutation relations
for ${\cal SL}_q$ in terms of three independent quantities $b$, $c$ and
$h$
\be
\begin{array}{l}
q^{2}hb - bh=2_q b \\
\rule{0pt}{5mm}
hc - q^{2}ch= -2_q c \\
\rule{0pt}{5mm}
\displaystyle
q(bc-cb)=h(1-\frac{q\lambda}{2_q}\,h).
\end{array}
\label{ex:sl-rea}
\ee
Note, that at $q\rightarrow 1$ relations (\ref{ex:sl-rea}) transforms
into the well  known
commutation relations for the generators of the universal enveloping
algebra $U(sl_2)$.

Starting from (\ref{ex:rep}) we have due to (\ref{rep-red})
$$
\bar\pi(h) = \xi\left(\matrix{q&0\cr
0&-q^{-1}}\right),\;\;
\bar\pi(b) = \xi\left(\matrix{0&q^{-1}\cr
0&0}\right),\;\;
\bar\pi(c) = \xi\left(\matrix{0&0\cr
q&0}\right),\quad \xi=\frac{q^2+1}{q^4+1}.
$$
At the classical limit this representation turns into the fundamental
vector representation of $U(sl_2)$. 

Let us point out that at $q\not=\pm 1$ the usual trace of $\bar\pi(h)$
does not 
equal to zero (this property is restored at the classical limit only).
However, the 
quantum trace of the matrix  $\bar\pi(h)$ is equal to zero
$$
Tr(C\cdot\bar \pi(h)) = 0.
$$
It should be emphasized that in the above relation the matrix 
$C$ is used to deform the trace of {\it  operators} of ${\cal SL}_q$
representation,
while in (\ref{q-sled}) the quantum trace is taken in ${\rm Mat}_n({\cal
SL}_q)$ (or in
${\rm Mat}_n(\Lq)$).

So, in order to retain the traceless property of ${\cal SL}_q$
representation in the space $V$ 
one has to modify the definition of the operator trace in ${\rm
End}(V)$: $Tr\rightarrow Tr_q$. 
This fact is no mere chance. The usual tensor category like that of
modules over $U(sl_n)$ is not 
suitable as the representation category for the algebra $\Lq$ (or ${\cal
SL}_q$). The natural 
representation category for the mentioned algebras is some {\it
quasitensor\footnote{The notion 
of the quasitensor category as a category of finite dimensional
representations of a quasitriangular 
Hopf algebra was introduced in \cite{Resh}.} Schur-Weyl category}. It is
the quantum trace 
(contrary to the usual one) which turns out to be a natural morphism
closely connected with the 
structure of the Schur-Weyl category. The detailed description of the
category, the role of the 
quantum trace in it and the connection with REA are considered in
\cite{GLS1} and \cite{GLS2}.

To complete the section, we clarify the connection of $B$ type
representation (\ref{ex:rep}) of 
mREA (\ref{ex:mrea}) with the result of \cite{Kul}. Applying the inverse
linear shift of generators 
(see (\ref{shift})) to representation (\ref{ex:rep}) one gets the
representation of REA with 
quadratic relations (\ref{ex:rea}). The matrices thus obtained are
equivalent (connected by a 
similarity transform) to the matrices of the two-dimensional
representation derived in \cite{Kul}. 
This representation is the lowest one in the series of non-commutative
finite dimensional 
representations of (\ref{ex:rea}). The higher dimensional modules of
this series are equivalent 
to the $q$-symmetrical tensor powers of the $B$ type modules. Their
construction is considered in 
the next section. So, in the particular case of $U_q(sl_2)$ $R$-matrix
our approach gives the known 
result of \cite{Kul}.

\section{Higher dimensional modules of $B$ type}
\label{th}

Let us turn to the problem of tensor product of mREA modules.
We are mainly interested in the following questions. Firstly, given the
fundamental 
module $V$ of $B$ type, how to define an mREA module structure in the
tensor power $V^{\otimes k}$?
Secondly, into which irreducible higher dimensional modules one can
decompose 
the module $V^{\otimes k}$? At last, how one 
can decompose the tensor product of arbitrary higher  dimensional
modules (not only fundamental 
ones) into the sum of irreducible components? Here we propose answers to
these questions.

\subsection{Tensor product of fundamental modules}

When trying to define the mREA module structure on the tensor product of
fundamental modules 
one finds a serious difficulty. The matter is that it is not known if
mREA (\ref{mREA}) (as well 
as REA (\ref{RE})) possesses the {\it bialgebra} structure. As a
consequence, in algebra (\ref{mREA}) 
one cannot define the coproduct operation.

To clarify the importance of this operation, consider the case of the
universal enveloping 
algebra ${\cal U} =  U(sl_n)$ and dwell upon the definition of the $\cal
U$-module structure on 
the tensor product of fundamental modules in the framework of Weyl
approach \cite{W}. 

The algebra $\cal U$ is a bialgebra\footnote{Moreover, the algebra $\cal
U$ (as the universal 
enveloping algebra of any Lie algebra) is a {\it Hopf} algebra. But for
our present purposes 
only the coproduct is needed.} with the cocommutative coproduct
$$
\Delta:\;\;{\cal U}\rightarrow {\cal U}\otimes {\cal U}.
$$
The action of $\Delta$ on a (Lie) generator\footnote{As is known, a Lie
algebra can be always 
embedded into its universal enveloping algebra. Here $x$ is the image of
a Lie algebra generator
under such an embedding.} $x$ of $\cal U$ is given by a simple formula
$$
\Delta(x) = x\otimes 1_{\cal U}+1_{\cal U}\otimes x\stackrel{\mbox{\tiny
def}}{=}
x_{(1)}+x_{(2)},
$$
where $1_{\cal U}$ is the unit element of $\cal U$.

Take now some irreducible representation $\rho:{\cal U}\rightarrow {\rm
End}(V)$ of ${\cal U}$ 
in a finite dimensional vector space $V$. To get a representation of
$\cal U$ in $V^{\otimes k}$ 
one first construct a homomorphism $\Delta^k: {\cal U}\rightarrow {\cal
U}^{\otimes k}$ by multiple
application of the coproduct $\Delta$. Then the image of $\cal U$ under
such a homomorphism is 
represented in  $V^{\otimes k}$ with the help of the map $\rho^{\otimes
k}$. For a 
given generator $x$ of $\cal U$ these two steps can be written in the
explicit form
\begin{eqnarray}
i)&& x\stackrel{\Delta^k}{\longrightarrow}\mathbf{x} =
x_{(1)}+x_{(2)}+\dots+x_{(k)}
\in {\cal U}^{\otimes k}\nonumber \\
\rule{0pt}{5mm}
ii)&&\mathbf{x}\rightarrow \rho^{\otimes k}(\mathbf{x})\in {\rm
End}(V^{\otimes k}).
\label{ro-k}
\end{eqnarray}

This representation is reducible. To extract the irreducible components
one uses the fact that 
in $V^{\otimes k}$ it is possible to define a natural representation of
the group algebra 
${\Bbb C}[{\cal S}_k]$ of the $k$-th order permutation group ${\cal
S}_k$. With respect to this
representation the space $V^{\otimes k}$ is decomposed into the direct
sum of irreducible 
${\Bbb C}[{\cal S}_k]$ modules $V_{\nu(a)}$. The modules are
parametrized by the standard Young 
tableaux $\nu(a)$, corresponding to all possible partitions $\nu\vdash
k$. With respect to
representation of $\cal U$ the subspaces $V_{\nu(a)}$ are also
irreducible. A generator
$x$ of $\cal U$ is represented in $V_{\nu(a)}$ by the following linear
operator
\be
\rho_{\nu(a)}(x) = P_{\nu(a)}\,\rho^{\otimes k}(\mathbf{x})\,P_{\nu(a)},
\label{nepr-class}
\ee
where $P_{\nu(a)}$ is the Young projector in $V^{\otimes k}$
corresponding to the tableau $\nu(a)$.
The modules parametrized by  different tableaux of the same partition
$\nu$ are
equivalent.

Return now to the case of mREA. As was already mentioned, this algebra
does not possess the coproduct
and for constructing the tensor product of mREA modules we cannot use
the above scheme as it stands. 
But it proves to be possible to generalize the final formulae
(\ref{ro-k}) and (\ref{nepr-class}) to 
the case of mREA.

Introduce the useful notation for a chain of $R$-matrices
\be
R_i\equiv R_{ii+1},\qquad
R^{\pm 1}_{(i\rightarrow j)}\stackrel{\mbox{\tiny def}}{=}
\left\{
\begin{array}{ll}
R^{\pm 1}_{i}R^{\pm 1}_{i+1}\dots R^{\pm 1}_{j}&\quad{\rm if}\;\;i<j\\
\rule{0pt}{5mm}
R^{\pm 1}_{i}R^{\pm 1}_{i-1}\dots R^{\pm 1}_{j}&\quad{\rm if}\;\;i>j\\
\rule{0pt}{5mm}
R^{\pm 1}_{i}&\quad {\rm if}\;\; i=j.
\end{array}
\right.
\label{R-chain}
\ee
The analog of reducible representation (\ref{ro-k}) is established in
the following proposition.

\begin{utv}\label{u-sec}
Consider the fundamental \Lq{} module V defined in Proposition \ref{u1}
and fix a basis $e_i$, 
$1\le i\le n$, in $V$. The tensor product $V^{\otimes k}$ is also an
\Lq{} module. In the basis 
$e_{i_1}\otimes \dots\otimes e_{i_k}$ of $V^{\otimes k}$ the matrices of
operators representing
the \Lq{} generators are  as follows
\be
\rho_k^t(l_i^{\,j}) = \pi^t(l_i^{\,j})\otimes I^{\otimes (k-1)}+ 
\sum_{s=1}^{k-1}R^{-1}_{(s\rightarrow 1)}\Bigl[\pi^t(l_i^{\,j})\otimes
I^{\otimes (k-1)}\Bigr]
R^{-1}_{(1\rightarrow s)}.
\label{ro-q}
\ee
Here $\rho_k^t$ and $\pi^t$ stand for the transposed matrices. 

\end{utv}

\noindent{\bf Proof\ \ }The proposition is proved by direct
calculations. Since the calculations are rather 
lengthy we shall not reproduce them in full detail, giving instead the
list of important intermediate steps 
with the corresponding results.

We prove the proposition by induction in $k$. For $k=1$ our assertion
reduces to Proposition \ref{u1}
and hence is true. Suppose it be true up to some $k-1$ and prove that
then it be true for $k$.

One should verify that operators (\ref{ro-q}) do satisfy the commutation
relations (\ref{mREA}).
Let us assign the number $k+1$ to the auxiliary space of indices of the
\Lq{} generators. One 
has to show 
\begin{eqnarray}
R_{k+1}\rho_k(L_{k+1})R_{k+1}\rho_k(L_{k+1}) -
\rho_k(L_{k+1})R_{k+1}\hspace*{-3.5mm}&&\hspace*{-3.5mm}
\rho_k(L_{k+1})R_{k+1} = \nonumber\\
&& R_{k+1}\rho_k(L_{k+1})-\rho_k(L_{k+1})R_{k+1}.
\label{doc-rep}
\end{eqnarray}
It is convenient to make the transposition of the matrices $\rho$ in the
above relation and to take into 
account that in accordance with (\ref{B-fun}) the matrix $\pi(l_i^{\,
j})$ reads
$$
\pi^t(L_{k+1})_1 = P_{1\,k+1}B_{k+1}, 
$$
where the unity enumerates the the matrix indices of the representation
space.
Then, on substituting (\ref{ro-q}) into (\ref{doc-rep}) we find that the
first summand 
in the left hand side decomposes into the sum of $k^2$ terms, a typical
one being as follows
$$
R_{k+1}R^{-1}_{(n\rightarrow 1)}P_{1\,k+1}B_{k+1}\Bigl[{\rm
Tr_{(1)}}R^{-1}_{(1\rightarrow n)}
R^{-1}_{(m\rightarrow 1)}B_1R_{1\,k+2} \Bigr]R^{-1}_{(1\rightarrow
m)}\equiv R_{k+1}Q(n,m),
$$
where the last equality is the definition of $Q(n,m)$, $0\le n,m\le
k-1$. Here in $Q(0,m)$ the 
unity matrices are substituted for the chains $R^{-1}_{(1\rightarrow
n)}$ and $R^{-1}_{(n\rightarrow 1)}$. 
The second summand in the left hand side of (\ref{doc-rep}) expands in a
similar  way but $R_{k+1}$ stands 
on the right of $Q$. 

In virtue of the supposition of the induction we conclude that it is
sufficient to
consider just $2k-1$ terms in each summand containing  $Q(k-1,n)$ and
$Q(n, k-1)$. So, one needs
to examine the following expression in the left hand side of
(\ref{doc-rep})
\begin{eqnarray}
\sum_{n=0}^{k-2}\Bigl(R_{k+1}(Q(n,k-1)+Q(k-1,n))\!\!\!&-&\!\!\!
(Q(n,k-1)+Q(k-1,n))R_{k+1}\Bigr)
\nonumber\\
&+&\!\!\!R_{k+1}Q(k-1,k-1)-Q(k-1,k-1)R_{k+1}.\label{st1}
\end{eqnarray}
The proposition will be proved if one could show that this is equal to
\be
R_{k+1}R^{-1}_{(k-1\rightarrow 1)}P_{1\,k+1}B_{k+1}R^{-1}_{(1\rightarrow
k-1)} -
R^{-1}_{(k-1\rightarrow 1)}P_{1\,k+1}B_{k+1}R^{-1}_{(1\rightarrow
k-1)}R_{k+1}.
\label{vklad}
\ee
Consider first the difference, containing $Q(n,k-1)$, $0\le n\le k-2$.
Since (see Appendix)
$$
{\rm Tr}_{(1)}R^{-1}_{(1\rightarrow n)}
R^{-1}_{(k-1\rightarrow 1)}B_1R_{1\,k+2} = R^{-1}_{(k-1\rightarrow
2)}P_{2\,k+2}B_{k+2}
R^{-1}_{(2\rightarrow n+1)}\quad \forall\,n\le k-2
$$ 
and (see, e.g., \cite{GPS})
\be
R_{12}B_1B_2=B_1B_2R_{12}
\label{RBB}
\ee
we find 
\begin{eqnarray}
R_{k+1}(Q(n,k-1)+
Q(k-1,n))\hspace*{-3mm}&-&\hspace*{-3mm}(Q(n,k-1)+Q(k-1,n))R_{k+1}=
\nonumber\\
\lambda R^{-1}_{(n\rightarrow 1)}R^{-1}_{(k-1\rightarrow 1)}
P_{1\,k+1}\hspace*{-3mm}&B&\hspace*{-4mm}_{k+1}P_{1\,k+2}B_{k+2}
R^{-1}_{(2\rightarrow k-1)} R^{-1}_{(1\rightarrow n)}
\label{Qnk}\\
&-&\hspace*{-3mm}\lambda R^{-1}_{(n\rightarrow
1)}R^{-1}_{(k-1\rightarrow 2)}
P_{1\,k+1}B_{k+1}P_{1\,k+2}B_{k+2} 
R^{-1}_{(1\rightarrow k-1)}R^{-1}_{(1\rightarrow n)}\nonumber
\end{eqnarray}
The calculation of the difference with $Q(k-1,k-1)$ in (\ref{st1}) is
more involved.
The trace contained in $Q(k-1,k-1)$ is as follows (see Appendix)
$$
{\rm Tr}_{(1)} R^{-1}_{(1\rightarrow k-1)}R^{-1}_{(k-1\rightarrow
1)}B_1R_{1\,k+2} = 
I -\lambda
P_{2\,k+2}B_{k+2}-\lambda\sum_{n=2}^{k-1}R^{-1}_{(n\rightarrow
2)}P_{2\,k+2}B_{k+2}
R^{-1}_{(2\rightarrow n)}.
$$
At last, it is a matter of straightforward calculation to show that the
unit matrix $I$ in the 
above expression leads to the necessary contribution (\ref{vklad}) while
the other terms exactly
cancel the unwanted summands of type
(\ref{Qnk}).\hfill\rule{6.5pt}{6.5pt}

\subsection{Decomposition of $V^{\otimes k}$ into mREA submodules}

Our next goal is to find the decomposition of the \Lq{} module
$V^{\otimes k}$ into irreducible 
submodules that is to find an analog of the classical formula
(\ref{nepr-class}). This can be 
done on the base of isomorphism $H_k(q)\cong {\Bbb C}[{\cal S}_k]$ which
was discussed in the first 
section. For our construction the most important consequence of such an
isomorphism is the existence 
of $q$-analogs of primitive Young idempotents ${\cal Y}_\nu(\sigma)\in
H_k(q)$ and decomposition
(\ref{razl-prostr}) of $V^{\otimes k}$ into the direct sum of $H_k(q)$
modules.

The main result of this section is formulated in the following
proposition.

\begin{utv}\label{u2}
Consider the mREA $\Lq$ generated by (\ref{mREA}) with the Hecke
$R$-matrix 
possessing the symmetry rank $p$ (see (\ref{Hec}), (\ref{closed}) and
(\ref{f-rank})).
Let  $V$ be the fundamental \Lq{} module of $B$ type with a fixed basis
$e_i$, 
$1\le i\le n$. According to Proposition \ref{u-sec} the space
$V^{\otimes k}$ is also 
an \Lq{} module for any positive integer $k$. Decompose the tensor
product $V^{\otimes k}$ 
into the direct sum (\ref{razl-prostr}). 

Then each component $V_{\nu(a)}$ of the direct sum  is an  \Lq{} module
and the generators 
$l_i^{\,j}$ are represented by linear operators
$\hat\pi_{\nu(a)}(l_i^{\,j})\in {\rm End}
(V_{\nu(a)})\hookrightarrow {\rm End}(V^{\otimes k})$. 
The matrices of these operators in the basis $e_{i_1}\otimes
\dots\otimes e_{i_k}$ of 
$V^{\otimes k}$ are of the form
\be
\pi^t_{\nu(a)}(l_i^{\,j}) =
Y_{\nu(a)}(R)\,\rho_k^t(l_i^{\,j})\,Y_{\nu(a)}(R),
\label{nepr-q}
\ee
where $\rho_k$ is defined in (\ref{ro-q}) of Proposition \ref{u-sec} and
the symbol $t$ means
the matrix transposition.

The modules parametrized by different tableaux of the same partition
$\nu\vdash k$
are equivalent.
\end{utv}

\noindent{\bf Proof\ \ }Consider matrices (\ref{ro-q}) of the \Lq{}
representation in $V^{\otimes k}$. 
Let us assign the number $k+1$ to the auxiliary space of indices of
\Lq{} generators and rewrite the 
commutation relations in form (\ref{doc-rep}). The projectors
$Y_{\nu(a)}(R)$ in (\ref{nepr-q}) are 
some polynomials in matrices $R_i$, $1\le i\le k-1$. The proof of the
proposition is based on the fact 
that the matrix $\rho^t_k(L_{k+1})$ commute with $R_i$ for $1\le i\le
k-1$. Indeed, taking this for 
granted for a moment, we conclude that $\rho^t_k(L_{k+1})$ commutes with
all projectors $Y_{\nu(a)}(R)$. 
Then, due to the orthogonality and the completeness of the set of
projectors
\begin{eqnarray}
&& Y_{\nu(a)}(R)\, Y_{\mu(b)}(R) =
\delta_{\nu\mu}\delta_{ab}Y_{\mu(b)}(R),\label{Y-norm}\\
&&\rule{0pt}{6mm}
\sum_{\nu\vdash k}\sum_{a=1}^{\dim[\nu]}Y_{\nu(a)}(R) = {\rm
id}_{V^{\otimes k}}\nonumber
\end{eqnarray}
it is easy to see that relation (\ref{doc-rep}) allows the projection
onto each component $V_{\nu(a)}$
of the direct sum (\ref{razl-prostr}) and the matrices of corresponding
representation are given 
by (\ref{nepr-q}).

To prove the commutativity supposed above one needs no particular
property of $\pi(l_i^{\,j})$ 
(an arbitrary $n\times n$ matrix $X$ can be substituted for
$\pi(l_i^{\,j})$). The proof 
consists in a simple calculation and completely based on the Yang-Baxter
equation (\ref{YBE}) 
and Hecke condition (\ref{Hec}).

The equivalence of modules parametrized by different Young tableaux of
the same partition $\nu$ 
stems from the connection of the corresponding Young idempotents. It can
be shown \cite{OgP} that 
any two idempotents ${\cal Y}_{\nu(a)}(\sigma)$ and ${\cal
Y}_{\nu(b)}(\sigma)$ are connected by 
a similarity transform. Namely, there exists an {\it invertible} element 
${\cal F}_{ab}(\nu|\sigma)\in H_k(q)$ which is a polynomial in
$\sigma_i$ such that
$$
{\cal F}_{ab}(\nu|\sigma)\,{\cal Y}_{\nu(a)}(\sigma)\,{\cal
F}_{ab}^{-1}(\nu|\sigma)  = 
{\cal Y}_{\nu(b)}(\sigma), \qquad 1\le a,b\le \dim[\nu].
$$
In the local representation (\ref{loc-rep}) of $H_k(q)$ in $V^{\otimes
k}$ the element 
${\cal F}_{ab}(\nu|\sigma)$ turns into the intertwining operator
$F_{ab}(R)$ which ensures the 
equivalence of $\pi_{\nu(a)}$ and $\pi_{\nu(b)}$
$$
F_{ab}\,\pi^t_{\nu(a)}\,F^{-1}_{ab} =
F_{ab}\,Y_{\nu(a)}\,\rho^t_k\,Y_{\nu(a)}\,F^{-1}_{ab} 
= Y_{\nu(b)}\, \rho^t_k\,Y_{\nu(b)} = \pi^t_{\nu(b)}.
$$
The second equality in the above chain of transformations is valid due
to commutativity
of $\rho^t_k$ and $F_{ab}(R)$ since the latter is a polynomial in $R_i$,
$1\le i\le k-1$.  
\hfill\rule{6.5pt}{6.5pt}
\medskip

A natural question arises about the irreducibility of modules
$V_{\nu(a)}$. Due to some reasons 
discussed at the end of the paper we suppose the modules to be
irreducible, but still we have 
no general proof of this fact. We state it as a quite plausible
hypothesis.

For the representations parametrized by single-row and single-column
diagrams the formulae
become much simpler and in this case one can explicitly calculate the
spectrum  of central 
elements (\ref{centr-el}). 

\begin{cor}\label{cor1}
Consider the \Lq modules $V_{\nu(a)}$, defined in Proposition \ref{u2}.
\begin{enumerate}
\item[{\bf i)}] For partitions $\nu=(k)$ and $\nu=(1^k)$ (for $k\le p$)
the matrices of operators representing the \Lq{} generators are given by
\begin{eqnarray}
&& \pi^t_{(k)}(l_i^{\,j}) = q^{1-k}k_q\,
S^{(k)}(R)\,\Bigl[\pi^t(l_i^{\,j})
\otimes I^{\otimes (k-1)}
\Bigr]\,S^{(k)}(R),\label{chi-s}\\
&&\rule{0pt}{5mm}
\pi^t_{[k]}(l_i^{\,j}) = q^{k-1}k_q\, A^{(k)}(R)\,\Bigl[\pi^t(l_i^{\,j})
\otimes I^{\otimes (k-1)}
\Bigr]\,A^{(k)}(R),\quad k\le p,
\label{chi-a}
\end{eqnarray}
where $S^{(k)}$ and $A^{(k)}$ are the $q$-symmetrizer and the
$q$-antisymmetrizer correspondingly.
\item[{\bf ii)}] In the representations $\pi_{(k)}$ and $\pi_{[k]}$ the
spectrum  $\chi$ of 
central elements $s_m = Tr_qL^m$ takes the following values
\begin{eqnarray}
&& \chi_{(k)}(s_m) = q^{-m(p+k-1)-p}\,k_q (p+k-1)_q^{m-1},
\label{char-s} \\
&& \rule{0pt}{5mm}
\chi_{[k]}(s_m) =q^{-m(p-k+1)-p}\,k_q (p-k+1)_q^{m-1}.
\label{char-a}
\end{eqnarray}
\end{enumerate}
\end{cor}

\noindent{\bf Proof\ \ }The assertion {\bf i)} is the direct consequence
of (\ref{nepr-q}). 
Indeed, taking into account the relations
$$
\begin{array}{l}
R_i^{\pm 1}S^{(k)}_{12\dots k} = q^{\pm 1}S^{(k)}_{12\dots k} = 
S^{(k)}_{12\dots k}R_i^{\pm 1}\\
\rule{0pt}{7mm}
R_i^{\pm 1}A^{(k)}_{12\dots k} = -q^{\mp 1}A^{(k)}_{12\dots k} = 
A^{(k)}_{12\dots k}R_i^{\pm 1}
\end{array}
\qquad 1\le i\le k-1
$$
and the definition (\ref{q-num}) of a $q$-number one immediately gets
(\ref{chi-s}) and 
(\ref{chi-a}) from (\ref{nepr-q}) where an arbitrary projector
$Y_{\nu(a)}$ should be 
replaced for $S^{(k)}(R)$ or $A^{(k)}(R)$ respectively. 

In order to find the values of characters (\ref{char-s}) and
(\ref{char-a}) we consider 
$\pi^t_{(k)}(Tr_qL^m)$  (take the $q$-symmetrical case for the
definiteness) and show that
this is a multiple of the $q$-symmetrizer (the identity operator on the
subspace $V_{(k)}$). The 
factor is equal to the character $\chi_{(k)}(s_m)$.

Taking into account (\ref{nonsing}) one rewrites the matrix involved in
the form 
$$
\pi^t_{(k)}(Tr_qL^m) = q^{m(1-k)-2p}\,k_q^m\,S^{(k)}_{12\dots
k}\Bigl[Tr_{(1)}B_1S^{(k)}_{12\dots k}
\Bigl]^{m-1}S^{(k)}_{12\dots k}.
$$
To calculate the trace $Tr_{(1)}B_1S^{(k)}_{12\dots k}$ we use the
recurrent
relations for the $q$-(anti)sym\-me\-tri\-zers (see, e.g., \cite{Gur})
\begin{eqnarray}
S^{(1)}(R) = I&\hspace*{5mm} &
S^{(k)}_{12\dots k}(R) = \frac{1}{k_q}\,S^{(k-1)}_{2\dots
k}(R)\Bigl(q^{1-k}I+
(k-1)_qR_{12}\Bigr)S^{(k-1)}_{2\dots k}(R),\label{S-pr}\\
A^{(1)}(R) = I&\hspace*{5mm} &
A^{(k)}_{12\dots k}(R) = \frac{1}{k_q}\,A^{(k-1)}_{2\dots
k}(R)\Bigl(q^{k-1}I-
(k-1)_qR_{12}\Bigr)A^{(k-1)}_{2\dots k}(R)\label{A-pr}
\end{eqnarray}
where $S^{(k-1)}_{2\dots k}$ is the $q$-symmetrizer of the $(k-1)$-th
order acting in components of 
$V^{\otimes k}$ with numbers from 2 till $k$. Then using (\ref{B-R})
and  (\ref{BC-norm}) one finds 
$$
Tr_{(1)}B_1S^{(k)}_{12\dots k} =
q^{-p}\,\frac{(p+k-1)_q}{k_q}\,S^{(k-1)}_{2\dots k}.
$$
Having found the above trace one gets  
\be
\pi^t_{(k)}(Tr_qL^m) =q^{-m(p+k-1)-p}\,k_q (p+k-1)_q^{m-1}\,S^{(k)}(R)
 = \chi_{(k)}(s_m)\,{\rm id}_{V_{(k)}},
\label{sym-cahr}
\ee
which proves (\ref{char-s}).  \hfill\rule{6.5pt}{6.5pt}
\medskip

It is worth pointing out that in the above formulae the central role
belongs to
the symmetry rank $p$ of the $R$-matrix but not to the
parameter\footnote{As was
already mentioned above there exist examples of $R$-matrices with
$p\not=n$
\cite{Gur}.} $n$ defining the size of $R$-matrix and the number of
generators in 
the corresponding REA \Lq. This feature is characteristic of the
Shur-Weyl category 
connected with the considered series of representations $\pi_{\nu}$
\cite{GLS1}.

\subsection{The sl-reduction}

As the next aim, we are going to explore a problem of tensor product of
fundamental 
${\cal SL}_q$ representations (\ref{rep-red}). Actually, the solution
for the problem 
follows from Proposition \ref{u2} and (\ref{sl-red-rep}). The only thing
we need is the 
spectrum of the element $Tr_qL$ in representations (\ref{nepr-q}).

\begin{lem} \label{l-chi-nu}
Let the partition $\nu\vdash k$ be of the height $s$ that is 
$$
\nu = (\nu_1,\nu_2,\dots,\nu_s),\quad\sum_{r=1}^{s}\nu_i = k, \quad
\nu_1\ge\nu_2\ge\dots\ge\nu_s>0.
$$
Then the spectrum of the central element $s_1 = Tr_qL$ in the
representation $\pi_{\nu(a)}$ 
$1\le a\le\dim[\nu]$ is as follows
\be
\chi_{\nu}(s_1) = q^{-2p}\sum_{r=1}^{s}q^{2r-1-\nu_r}(\nu_r)_q,
\label{nu-tr1}
\ee
where $p$ is the symmetry rank of $R$-matrix and $(\nu_r)_q$ is the
$q$-analog of the integer $\nu_r$
(see definition (\ref{q-num})). 
\end{lem}

\noindent{\bf Proof\ \ }Let us find the matrix $\pi_{\nu(a)}(Tr_qL)$ for
an arbitrary $\nu(a)$. 
Using (\ref{q-sled}) and  (\ref{nonsing}) one immediately gets from
(\ref{nepr-q})
$$
\pi^t_{\nu(a)}(Tr_qL) =q^{-2p}\,
Y_{\nu(a)}(R)\Bigl(J_1^{-1}+J_2^{-1}+\dots +
J_k^{-1}\Bigr)Y_{\nu(a)}(R).
$$
Here the matrices $J_k$ read
\be
J_1 = I,\quad J_i = R_{(i-1\rightarrow 1)}R_{(1\rightarrow i-1)}\quad
i\ge 2.
\label{J}
\ee
They are the images of the {\it Jucys-Murphy} elements ${\cal
J}_i(\sigma)$ of the Hecke algebra
$H_k(q)$ under the local representation (\ref{loc-rep}). The elements
${\cal J}_i(\sigma)$ generate the
maximal commutative subalgebra in $H_k(q)$ and have the following
important property (see, e.g., 
\cite{OgP})
$$
{\cal J}_i(\sigma){\cal Y}_{\nu(a)}(\sigma) = {\cal
Y}_{\nu(a)}(\sigma){\cal J}_i(\sigma)
 = q^{2(c_i - r_i)}{\cal Y}_{\nu(a)}(\sigma),
$$
where $c_i$ and $r_i$ are respectively the coordinates of the column and
the row to which
the box with the number $i$ belongs. The quantity $q^{2(c-r)}$ is called
a {\it content} of the
$(c,r)$-th box of a given Young diagram. Here is an example of the
diagram $\nu=(4,3,1^2)$
with corresponding contents
$$
\begin{array}{|c|c|c|c|}\hline
1&q^2&q^4&q^6\\ \hline
q^{-2}&1&q^2&\multicolumn{1}{c}{}\\ \cline{1-3}
q^{-4}&\multicolumn{3}{c}{}\\ \cline{1-1}
q^{-6}&\multicolumn{3}{c}{}\\ \cline{1-1}
\end{array}
$$

Therefore we come to the result 
\be
\pi^t_{\nu(a)}(Tr_qL) = q^{-2p}
\Bigl(\sum_{i=1}^{k}q^{-2(c_i-r_i)}\Bigr)Y_{\nu(a)} = \chi_{\nu}(Tr_qL)
Y_{\nu(a)}.
\label{chi1-nu}
\ee
It is obvious that the sum of all contents (or their inverses) of a
given tableau
$\nu(a)$ is completely defined by the partition (diagram) $\nu$ and has
the same value for 
all tableaux $\nu(a)$ corresponding to a given $\nu$. As a consequence,
the character
$\chi_{\nu}(Tr_qL)$ does not depend on $a$. Of course, this is in
agreement with the 
equivalence of representations $\pi_{\nu(a)}$ parametrized by different
tableaux $\nu(a)$ of 
the same partition $\nu$.

It is a matter of a simple calculation to show that $\chi_{\nu}(Tr_qL)$
represented as the sum of all inverse contents of the tableau 
$\nu_{(a)}$ can  be written in form
(\ref{nu-tr1}).\hfill\rule{6.5pt}{6.5pt}
\medskip

It is interesting to point out one peculiarity of the quantum case.
Namely, the spectrum of
central elements distinguishes the representations more effectively than
in the classical
case. To show this, take the $q$ as a parameter and extend the complex
field
$\Bbb C$ to the field of rational functions in $q$. The character
$\chi_\nu(s_1)$ can be
identically transformed to the following expression
$$
q^{2p}\lambda\,\chi_{\nu}(s_1) = q^ss_q-\sum_{r=1}^{s}q^{2(r-\nu_r)-1}
\equiv q^ss_q - \Omega(\nu_1,\dots,\nu_s).
$$
Given $\Omega(\nu_1,\dots, \nu_s)$ as a function in $q$, one
unambiguously restores
the values of all $\nu_r$, since the numbers $2(r-\nu_r)$ form a {\it
strictly} 
increasing sequence. Therefore if the height $s$ of a partition $\nu$ is
fixed and
known then the spectrum of the first central element $s_1 = Tr_q L$ is
sufficient
to distinguish the representations. It is not so, of course, for the
classical case
$q\rightarrow 1$.

Now we are in a position to formulate the results for the higher
dimensional modules of the 
algebra ${\cal SL}_q$.

\begin{utv}\label{u3}
Consider the algebra ${\cal SL}_q$ (\ref{sl-red}) with a Hecke
$R$-matrix possessing the 
symmetry rank $p$. Let  $V$ be the fundamental \Lq{} module of $B$ type
with a fixed basis 
$e_i$, $1\le i\le n$. According to Proposition \ref{u-sec} the tensor
product $V^{\otimes k}$ 
is also an \Lq{} module for any $k\in {\Bbb N}$.  
The following assertions are true.

\begin{enumerate}
\item[{\bf i)}] The space $V^{\otimes k}$ is a reducible ${\cal SL}_q$
module.
In the basis $e_{i_1}\otimes\dots\otimes e_{i_k}$ the matrices of
operators 
representing the ${\cal SL}_q$ generators $f_i^{\, j}$ have the
following form
\be
\bar\rho^t_k(f_i^{\, j}) = \frac{1}{\omega}
\Bigl(\rho^t_{k}(l_i^{\,j}) - \frac{\delta_i^{\,j}}{p_q q^p}{\cal
Z}_k\Bigr),\qquad
\omega = \frac{q^{1-p}}{p_q}\,(q^{p-2}(p+1)_q - 1),
\label{sl-ro-q}
\ee
where the symbol $t$ means the matrix transposition, $\rho_k(l_i^{\,j})$
is defined in (\ref{ro-q}) 
and ${\cal Z}_k$ is given by
$$
{\cal Z}_k = I +\sum_{n=1}^{k-1}R^{-1}_{(n\rightarrow
1)}R^{-1}_{(1\rightarrow n)}.
$$

\item[{\bf ii)}]
Decompose the tensor product $V^{\otimes k}$ into the direct sum
(\ref{razl-prostr}).
Each  component $V_{\nu(a)}$ of the direct sum is an ${\cal SL}_q$
submodule in 
$V^{\otimes k}$. The generators $f_i^{\,j}$ are represented by linear
operators with the 
following matrices
\be 
\bar\pi^t_{\nu(a)}(f_i^{\,j}) = \frac{1}{\omega_\nu}\,
Y_{\nu(a)}\Bigl[\rho_k^t(l_i^{\,j}) -
\delta_i^{\,j}\,\frac{q^p}{p_q}\,\chi_\nu(s_1)
I_{12\dots k}\Bigr]Y_{\nu(a)},\qquad \omega_\nu =
1-\lambda\,\frac{q^p}{p_q}\,\chi_\nu(s_1),
\label{sl-nepr-q}
\ee
the character $\chi_\nu(s_1)$ being defined in (\ref{nu-tr1}).
The modules parametrized by different tableaux $\nu(a)$ of the same
partition $\nu\vdash k$
are equivalent.

\item[{\bf iii)}] The spectrum $\bar\chi$ of the ${\cal SL}_q$ central
elements 
$\bar s_m = Tr_qF^m$ in representations $\bar\pi_{(k)}$ and
$\bar\pi_{[k]}$ corresponding
to $\nu = (k)$ and $\nu=(1^k)$ takes the following values
\be
\bar\chi_{(k)}(\bar s_m) =
q^{-p-m}\,\frac{k_q(p-1)_q(p+k)_q}{(p+k-1)_q}\,
\frac{\Bigl[(p-1)_q^{m-1}(p+k)_q^{m-1}+(-1)^mk_q^{m-1}\Bigr]}{(q^{p-2}(p+k)_q-k_
q)^m}
\label{slchi-s}
\ee
\be
\bar\chi_{[k]}(\bar s_m) =
q^{-p+m}\,\frac{k_q(p+1)_q(p-k)_q}{(p-k+1)_q}\,
\frac{\Bigl[(p+1)_q^{m-1}(p-k)_q^{m-1}+(-1)^mk_q^{m-1}\Bigr]}{(q^{p+2}(p-k)_q+k_
q)^m}
\label{slchi-a}
\ee
\end{enumerate}
\end{utv}

\noindent{\bf Proof\ \ }This proposition is the direct consequence of
Propositions \ref{u-sec}
and \ref{u2}, Corollary \ref{cor1} and rule (\ref{sl-red-rep}). Indeed,
the operators 
$\bar\rho_k(Tr_qF)$ and $\bar\pi_{\nu(a)}(Tr_qF)$ are obviously equal to
zero. Basing on 
Propositions \ref{u-sec} and  \ref{u2} one can verify that
(\ref{sl-ro-q}) and (\ref{sl-nepr-q}) 
do satisfy (\ref{mREA}) and that the factor $\omega_\nu^{-1}$ ensures
the proper normalization 
of the right hand side of (\ref{mREA}). 

As for the values (\ref{slchi-s}) and (\ref{slchi-a}), they can be found
by straightforward 
but rather lengthy calculations on the base of (\ref{sl-nepr-q}),
(\ref{char-s}) and 
(\ref{char-a}).
\hfill \rule{6.5pt}{6.5pt}
\medskip

In the case of $U_q(sl_n)$ $R$-matrix we have $p=n$ and at the classical
limit  $q\rightarrow 1$ 
the spectrum (\ref{slchi-s}), (\ref{slchi-a}) of the ${\cal SL}_q$
central elements tends to 
the spectrum of the $U(sl_n)$ Casimir elements in the corresponding
representations (see, e.g.,
\cite{BRon}).

At last, consider the tensor product of two (irreducible) modules
$V_\mu$ and $V_\nu$ over 
$\Lq$ (or ${\cal SL}_q$) in which the representation operators $\pi_\mu$
and $\pi_\nu$ are
given by (\ref{nepr-q}) (or by (\ref{sl-nepr-q})). Using the isomorphism
$H_k\cong 
{\Bbb C}[{\cal S}_k]$ one can show (see, e.g., \cite{GLS1}) that the
tensor product 
$V_\mu\otimes V_\nu$ is isomorphic to the following direct sum of \Lq{}
(or ${\cal SL}_q$)
modules $V_\sigma$
\be
V_\mu\otimes V_\nu \cong c_{\mu\nu}^{\sigma}V_\sigma.
\ee
Here $c_{\mu\nu}^{\sigma}$ are the Littlewood-Richardson coefficients
defining a ring structure 
in the set of Schur symmetric functions.

\section{Fundamental module of $R$ type }
\label{f}

Besides the fundamental REA module of $B$ type considered in the
previous sections, one can 
construct another module which will be called {\it the fundamental
module of $R$ type}. In 
the case of $R$-matrix connected with $U_q(sl_n)$ the corresponding
representation originates 
from the general theory of dual Hopf algebras and we generalize it to
the case of an 
arbitrary $R$-matrix.

\subsection{The definition and tensor product decomposition rule}

So, suppose at first that $R$-matrix defining the structure of REA is
the
Drinfeld-Jimbo $R$-matrix connected with the quantum universal
enveloping algebra $U_q(sl_n)$.
In this case the Hopf algebra (\ref{RTT}) is an algebra ${\rm
Fun}_q(GL(n))$ of functions on 
the quantum group \cite{FRT}. Besides, there exists embedding
(\ref{mapinq}) of the corresponding
REA $\Lq$ into $U_q(gl_n)$ --- the dual Hopf algebra to ${\rm
Fun}_q(GL(n))$. As a consequence,
it is possible to define a pairing among the generators $\hat l_i^{\,j}$
and $t_i^{\,j}$. 

Using the explicit formulae for the paring of $U_q(gl_n)$ generators
$L^\pm$ and
${\rm Fun}_q(GL(n))$ generators $T$ (see \cite{FRT}), we get the
following result
\be
{\cal h} T_1T_2\dots T_k , \hat L_{k+1} {\cal i} = R_{(k\rightarrow
1)}R_{(1\rightarrow k)}
\equiv J_{k+1}.
\label{LT-par}
\ee
Here we have used the compact notations (\ref{R-chain}) and (\ref{J})
for the chains of $R$-matrices.

As is known from the Hopf algebra theory, any module over a Hopf algebra
$\cal H$ can be transformed
into a comodule over its dual Hopf algebra ${\cal H}^*$ and vice versa.
Let us use this fact in order
to define a representation of REA \Lq{} in a finite dimensional vector
space $V$, $\dim V =n$.

On fixing a basis $e_i$, $1\le i\le n$, one can convert the space $V$
into a {\it left} comodule
over the Hopf algebra ${\cal H} = {\rm Fun}_q(GL(n))$ by means of the
corepresentation $\delta$
$$
\delta:\;V\rightarrow {\cal H}\otimes V,\qquad \delta(e_i) =
t_i^{\,j}\otimes e_j
$$
where the summation over the repeated indices is understood.
Then in $V$ we get a {\it right} \Lq{} action by the following rule
$$
e_1\triangleleft \hat L_2 = {\cal h}T_1,\hat L_2{\cal i}\, e_1 = 
R_{12}^2\,e_1.
$$

This action can be easily expanded to the tensor product $V^{\otimes k}$
for 
any $k\in {\Bbb N}$. Indeed the comodule structure of $V^{\otimes k}$ is
obvious
$$
\delta_k:\; V^{\otimes k}\rightarrow {\cal H}\otimes V^{\otimes
k},\qquad
e_1\otimes \dots \otimes e_k\stackrel{\delta_k}{\longrightarrow}T_1\dots
T_k\otimes 
(e_1\otimes \dots \otimes e_k).
$$
Hence, taking into account (\ref{LT-par})
$$
e_1\otimes \dots \otimes e_k\triangleleft \hat L_{k+1} = 
J_{k+1} \,e_1\otimes \dots \otimes e_k.
$$

It turns out that we can directly generalize the above formulae to the
case of
an {\it arbitrary} $R$-matrix. The following proposition is easy to
verify.

\begin{utv}
\label{u4}
Consider the REA $\Lq$ generated by relations (\ref{RE}) with an
arbitrary $R$-matrix. 
The matrix will be treated as that of a linear operator acting in the
tensor square of 
a finite dimensional vector space $V$, $\dim V = n$. Define a linear map
$\theta_k: 
\Lq\rightarrow {\rm End}(V^{\otimes k})$ by the following rule 
\be\left\{\begin{array}{l}
\theta_k(e_{\cal L}) = {\rm id}_{V^{\otimes k}}\\
\rule{0pt}{5mm}
\theta_k(\hat L_{k+1}) = \alpha\,J_{k+1}\\
\rule{0pt}{5mm}
\theta_k(\hat l_1\cdot \hat l_2\cdot\dots\cdot \hat l_m) = \theta_k(\hat
l_1)\cdot 
\theta_k(\hat l_2)\cdot\dots\cdot \theta_k(\hat l_m),
\end{array}
\right.
\label{R-rep}
\ee
where $\alpha\not=0$ is an arbitrary complex number. Then $\theta_k$
realizes a representation
of $\Lq$ in the space $V^{\otimes k}$.
\end{utv}

\noindent{\bf  Proof}\ \  It is sufficient to substitute the matrices
$\theta_k(\hat L_{k+1}) 
= \alpha J_{k+1}$ in (\ref{RE}) rewritten in the form
$$
R_{k+1}\hat L_{k+1}R_{k+1}\hat L_{k+1} - \hat L_{k+1}R_{k+1}\hat
L_{k+1}R_{k+1} = 0
$$
and make use of the following consequence of Yang-Baxter equation
(\ref{YBE}) 
\be
\begin{array}{l}
(R_1\dots R_k)R_i = R_{i+1}(R_1\dots R_k), \quad 1\le i\le k-1\\
\rule{0pt}{5mm}
(R_k\dots R_1)R_i = R_{i-1}(R_k\dots R_1), \quad 2\le i\le k.
\end{array}
\label{R-ch}
\ee
As for the numeric factor $\alpha\not=0$, it can be arbitrary due to the
renormalization
automorphism $\hat L\rightarrow\alpha \hat L$ (see Remark \ref{r1}).
\hfill\rule{6.5pt}{6.5pt}
\medskip

Note that in proving Proposition \ref{u4} we use nothing but the
Yang-Baxter equation 
for the $R$-matrix. Therefore, representation (\ref{R-rep}) is valid not
only for the
quantum group $R$-matrix but also for an arbitrary solution of the
Yang-Baxter 
equation (even of a non-Hecke type). 

At $k>1$ the representation $\theta_k$ is reducible. The Hecke condition
(\ref{Hec}) 
is needed for extracting the irreducible components of $\theta_k$.

\begin{utv}\label{u5}
Let REA $\Lq$ be generated by (\ref{RE}) with a Hecke $R$-matrix.
Consider the representation $\theta_k$ (\ref{R-rep}) in the space
$V^{\otimes k}$. 
Decompose $V^{\otimes k}$ into the direct sum of subspaces $V_{\nu(a)}$
in accordance with
(\ref{razl-prostr}). 

Then each $V_{\nu(a)}$ is an \Lq{} submodule and the matrices of linear
operators representing
the generators $\hat l_i^{\,j}$ are given by 
\be
\theta_{\nu(a)}(\hat L_{k+1}) = Y_{\nu(a)}(R)\,\theta_k(\hat L_{k+1})\,
Y_{\nu(a)}(R).
\label{th-lam}
\ee
The modules parametrized by different tableaux of the same partition
$\nu\vdash k$ 
are equivalent.
\end{utv}

\noindent{\bf Proof\ \ } The proof is based on the fact that
$q$-projectors $Y_\nu(R)$
are actually polynomials in $J_1,\dots ,J_{k}$ for $\nu\vdash k$ (see
\cite{OgP}).
Being the images of Jucys-Murphy elements, the operators $J_i$
commute\footnote{
One can verify this fact independently, on the base of (\ref{R-ch}).}
with $J_{k+1}$. As a consequence, the relation 
$$
R_{k+1}\theta_k(\hat L_{k+1})R_{k+1}\theta_k(\hat L_{k+1}) =
\theta_k(\hat L_{k+1})R_{k+1}\theta_k(\hat L_{k+1})R_{k+1} 
$$
which takes place in  ${\rm End}(V^{\otimes k})$ admits projection into
subspaces 
${\rm End}(V_\nu)$ in accordance with (\ref{th-lam}).

The equivalence of $V_{\nu(a)}$ and $V_{\nu(b)}$ corresponding to
different tableaux
of the same partition $\nu$ is proved in the same way as in Proposition
\ref{u2}.
\hfill\rule{6.5pt}{6.5pt}
\medskip

It is worth mentioning, that for constructing representations
$\theta_{\nu}$ one has no need
the symmetry rank of $R$ to be finite. 

As in the case of $B$ type module one can explicitly calculate the
spectrum of
central elements (\ref{centr-el}) in the representations parametrized by
single-row and single-column diagrams.

\begin{cor} In the representations parametrized by partitions $\nu=(k)$
and
$\nu=(1^k)$ the spectrum $\hat \chi$ of the central elements $s_m =
Tr_q\hat L^m$
takes the following values
\begin{eqnarray}
&&\hat\chi_{(k)}(s_m) =
q^{-p}\Bigl(q^{-2m}p_q+\lambda\,\frac{(p+k)_q}{(k+1)_q}\,
q^{m(k-1)}[m(k+1)]_q\Bigr)\label{chi-th-s}\\
&&\hat\chi_{[k]}(s_m) =
q^{-p}\Bigl(q^{2m}p_q-\lambda\,\frac{(p-k)_q}{(k+1)_q}\,
q^{-m(k-1)}[m(k+1)]_q\Bigr)\qquad k\le p
\end{eqnarray}
\end{cor}

\noindent{\bf Proof\ \ } We shall consider the case of $q$-symmetric
representation
$\theta_{(k)}$ for the definiteness. Let us first calculate
$\theta_{(k)}(Tr_q\hat L)$.
In accordance with (\ref{R-rep}) and (\ref{th-lam}) the matrices
representing
REA generators read
$$
\theta_{(k)}(\hat L_{k+1}) = S^{(k)}J_{k+1}S^{(k)}.
$$
Taking into account (see \cite{OgP}) that
\be
S^{(k+1)} = S^{(k)}\,\frac{J_{k+1}-q^{-2}}{q^{2k} - q^{-2}}
\label{sym-J}
\ee
we rewrite the matrix $\theta_{(k)}(\hat L_{k+1})$ in the equivalent
form
$$
\theta_{(k)}(\hat L_{k+1}) = \lambda q^{k-1}(k+1)_q
S^{(k+1)}+q^{-2}S^{(k)}.
$$
Next, in virtue of
$$
{Tr_q}_{(k+1)}S^{(k+1)} = q^{-p}\,\frac{(p+k)_q}{(k+1)_q}\,S^{(k)}
$$
we come to the final result
\be
\theta_{(k)}(Tr_q \hat L) = q^{-p}\left(q^{-2}p_q+\lambda
q^{k-1}(p+k)_q\right)
\,S^{(k)}\equiv \hat \chi_{(k)}(s_1)\,{\rm id}_{V_{(k)}}.
\label{trace-1}
\ee
Then, basing on (\ref{sym-J}) and on the relations 
$$
S^{(k+1)}S^{(k)} = S^{(k+1)},\qquad S^{(k+1)}J_{k+1} = q^{2k}S^{(k+1)}
$$
one can prove (\ref{chi-th-s}) by induction in the power 
$m$ of $Tr_q\hat L^m$, where (\ref{trace-1}) serves as 
the first step. The final step of induction gives
$$
\hat \chi_{(k)}(s_m) = q^{-p}\Bigl(q^{-2m}p_q+\lambda q^{m(k-1)}(p+k)_q
\sum_{r=0}^{m-1}q^{(k+1)(2r+1-m)}\Bigr).
$$
With substitution $t=q^{k+1}$ one can easily show that
$$
\sum_{r=0}^{m-1}q^{(k+1)(2r+1-m)} = \frac{[m(k+1)]_q}{(k+1)_q}
$$
coming thereby to the desired result
(\ref{chi-th-s}).\hfill\rule{6.5pt}{6.5pt}
\medskip

The representation  $\theta_1$ in the space $V$ itself is irreducible
and the matrices
of operators representing \Lq{} generators are as follows
\be
\theta_1(\hat L_2) = R^2_{12}.
\label{r2-fun}
\ee
where the indices of the first space are those of matrices from ${\rm
End}(V)$ and 
the indices of the second space enumerates the generators of the
algebra. The space 
$V$ with the above \Lq{} representation will be called {\it the
fundamental module 
of $R$ type}.

The representation of mREA (\ref{mREA}) obtained from (\ref{r2-fun}) by
shift (\ref{shift})
reads
\be
\theta_1(L_2) = - R_{12},
\label{sl-r2}
\ee
where we first perform a renormalization of (\ref{r2-fun}) by the factor 
$\alpha = -\lambda^{-1}$.

\subsection{The sl-reduction}

In order to pass from the REA representation $\theta_{\nu(a)}$
(\ref{th-lam}) to the corresponding
representation $\bar\theta_{\nu(a)}$ of the algebra ${\cal SL}_q$ we
need to calculate the
spectrum of $\theta_{\nu(a)}(Tr_q\hat L)$. 

\begin{lem}
Let the partition $\nu\vdash k$ be of the height $s$ that is 
$$
\nu = (\nu_1,\nu_2,\dots,\nu_s),\quad\sum_{r=1}^{s}\nu_i = k, \quad
\nu_1\ge\nu_2\ge\dots\ge\nu_s>0.
$$
Then the spectrum of the central element $s_1 = Tr_q\hat L$ in the
representation $\theta_{\nu(a)}$ 
$1\le a\le\dim[\nu]$ is as follows
\be
\theta_{\nu(a)}(Tr_q\hat L) = \zeta_\nu(s_1)Y_{\nu(a)},\qquad
\zeta_\nu(s_1) = q^{-p}p_q+\lambda\sum_{r=1}^{s}q^{\nu_r+1-2r}(\nu_r)_q,
\ee
where $p$ is the symmetry rank of $R$-matrix and $(\nu_r)_q$ is the
$q$-analog of the integer $\nu_r$
(see definition (\ref{q-num})). 
\end{lem}

\noindent{\bf Proof\ \ } The lemma is proved by the straightforward
calculation in analogy with
the proof of Lemma \ref{l-chi-nu}.\hfill\rule{6.5pt}{6.5pt}

Let us consider the representations $\theta_{(k)}$ and $\theta_{[k]}$
corresponding to 
single-row and single-column diagrams in more detail. In this case one
can explicitly calculate 
the spectrum of central elements similarly to the $B$ type
representation.

\begin{utv}
Let the Hecke type $R$-matrix has the symmetry rank $p$. Consider the
REA representations of $R$ type
$\theta_{(k)}$ and $\theta_{[k]}$ parametrized by partitions $\nu=(k)$
and $\nu=(1^k)$ (\ref{th-lam}).
Then the corresponding representations $\bar\theta$ of the ${\cal SL}_q$
generators $f_i^{\,j}$
are as follows
\begin{eqnarray}
&&
\bar\theta_{(k)}(F_{k+1}) = \frac{q^{1-p}(p+k)_q}{q^{2-p}(p+k)_q -
k_q}\,\Bigl(S^{(k)}I_{k+1} - 
\frac{p_q(k+1)_q}{(p+k)_q}\,S^{(k+1)}\Bigr)\\
&&
\bar\theta_{[k]}(F_{k+1}) = \frac{q^{-1-p}(p-k)_q}{q^{-2-p}(p-k)_q +
k_q}\,\Bigl(
\frac{p_q(k+1)_q}{(p-k)_q}\,A^{(k+1)} - A^{(k)}I_{k+1}\Bigr)
\end{eqnarray}
The spectrum $\bar\zeta$ of the ${\cal SL}_q$ central elements $\bar s_m
= Tr_qF^m$ in these 
representations takes the following values
\begin{eqnarray}
&&
\bar\zeta_{(k)}(\bar s_m) =
q^{-p-m(p-1)}\,\frac{k_q(p-1)_q(p+k)_q}{(k+1)_q}\,
\frac{\Bigl[(p+k)_q^{m-1}+(-1)^{m}k_q^{m-1}(p-1)_q^{m-1}\Bigr]}{(q^{2-p}(p+k)_q
- k_q)^m}\\
&&
\bar\zeta_{[k]}(\bar s_m) =
q^{-p-m(p+1)}\,\frac{k_q(p+1)_q(p-k)_q}{(k+1)_q}\,
\frac{\Bigl[(-1)^m(p-k)_q^{m-1}+k_q^{m-1}(p+1)_q^{m-1}\Bigr]}{(q^{-2-p}(p-k)_q
+ k_q)^m}
\end{eqnarray}
\end{utv}

\noindent{\bf Proof\ \ } The proof consists in direct calculations on
the base of 
(\ref{sl-red-rep}) and we shall not present it
here.\hfill\rule{6.5pt}{6.5pt}

\subsection{Interrelation between modules of $B$ and $R$ types}

Let us now find a connection between the fundamental modules of $B$ and
$R$ types.
If the symmetry rank $p=2$ (for example, when $R$ stems from
$U_q(sl_2)$)
these modules are equivalent. To be more precise, the situation is as
follows.
In virtue of (\ref{f-rank}) the $q$-antisymmetrizer $A^{(2)}$ is a unit
rank
projector in $V^{\otimes 2}$ and its matrix can be written in the form
$$
A_{i_1i_2}^{\;j_1j_2} = u_{i_1i_2}v^{j_1j_2},
$$
the matrices $\|u_{ij}\|$ and $\|v^{ij}\|$ being nonsingular. Then one
can show 
that the representations $\pi$ and $\theta_1$ of mREA (\ref{mREA}) are
connected by 
the relation
$$
q^2\,u_1\cdot \pi(L_2)\cdot u_1^{-1} = qI_{12} + \theta_1(L_2).
$$
After the sl-reduction we come to the representations $\bar\pi$ and
$\bar\theta_1$ 
of the algebra ${\cal SL}_q$ (\ref{sl-red}) and simplify the above
formula to the 
expression
$$
u_1\cdot\bar \pi(F_2)\cdot u_1^{-1} = \bar\theta_1(F_2),
$$
$F=\|f_i^{\,j}\|$ being the matrix composed of the ${\cal SL}_q$
generators.

In the case $p>2$ the fundamental modules of $B$ and $R$ types are not
equivalent. Constraining 
ourselves to the case of ${\cal SL}_q$ algebra (\ref{sl-red}) we shall
prove that $R$ type
representation $\bar \theta(F)$ is equivalent to $\bar\pi_{[p-1]}(F)$
obtained from (\ref{chi-a})
by means of sl-reduction (\ref{sl-nepr-q}).

For this purpose, consider in more detail the structure of the subspace
$V_{[p-1]}\subset 
V^{\otimes (p-1)}$. By definition (\ref{razl-prostr}) the subspace
$V_{[p-1]}$ is the image 
of the $q$-antisymmetrizer $A^{(p-1)}$
$$
V_{[p-1]} = A^{(p-1)}(R)\triangleright V^{\otimes (p-1)}.
$$ 
Since the symmetry rank of $R$-matrix is equal to $p$ then the
$q$-antisymmetrizer $A^{(p)}$ 
is a unit rank projector and  its matrix can be written in the form
\be
{A^{(p)}}_{i_1\dots i_p}^{\;j_1\dots j_p} = u_{i_1\dots i_p}v^{j_1\dots
j_p}
\label{A-uv}
\ee
where as follows from (\ref{Y-norm}) the tensors $u$ and $v$ are
normalized by the condition 
$$
\sum_{\{i\}}u_{i_1\dots i_p}v^{i_1\dots i_p} = 1.
$$
It is convenient to introduce the following linear combinations of the
basis vectors of
the space $V^{\otimes (p-1)}$
\be
\epsilon^i\stackrel{\mbox{\tiny def}}{=}\sum_{\{a\}}v^{ia_2\dots
a_p}e_{a_2}\otimes\dots\otimes e_{a_p}.
\label{dual-bas}
\ee
The following lemma establishes an important property of the vectors
$\epsilon^i$.

\begin{lem} \label{l3}
Consider the set of $n$ vectors $\epsilon^i\in V^{\otimes (p-1)}$
defined in (\ref{dual-bas}).
These are eigenvectors of the $q$-antisymmetrizer $A^{(p-1)}$ and they
form a basis
of the subspace $V_{[p-1]}$
\be
A^{(p-1)}(R)\triangleright \epsilon^i = \epsilon^i,\qquad  
\forall\,\mathbf{w}\in V_{[p-1]}:\quad \mathbf{w} = \sum_i
w_i\epsilon^i.
\ee
\end{lem}

\noindent{\bf Proof\ \ } Consider the recurrence relation (\ref{A-pr})
for the 
$q$-antisymmetrizer $A^{(p)}$ and calculate the trace in the first
matrix space 
with the matrix $B_1$
$$
Tr_{(1)}B_1A^{(p)}_{12\dots p} = \frac{1}{q^pp_q}\,A^{(p-1)}_{2\dots p}.
$$
In virtue of (\ref{A-uv}) we get the following expression for the matrix
$\Bbb A$ 
of $A^{(p-1)}$
\be
{\Bbb A}_{a_2\dots a_p}^{\;\;b_2\dots b_p} =q^pp_q
\sum_{m,n}B_m^{\,n}u_{na_2\dots a_p}
v^{mb_2\dots b_p},
\label{*}
\ee
where for the compactness we omit the superscript $(p-1)$ of the matrix
$\Bbb A$.
Also we need the formula connecting the matrix $C$ (\ref{BC}) and the
tensors $u$ 
and $v$. It can be shown (see \cite{Gur}) that
\be
C_i^{\,j} = \frac{p_q}{q^p}\,\sum_{\{a\}}u_{ia_2\dots a_p}v^{ja_2\dots
a_p}\equiv
\frac{p_q}{q^p}\,\sum_{\{a\}}u_{i\{a\}}v^{j\{a\}},
\label{C-uv}
\ee
where in the last equality we have introduced a convenient multi-index
notation.

At last, taking into account definition (\ref{dual-bas}) we get the
necessary result 
(the summation  over the repeated indices is understood)
\begin{eqnarray*}
A^{(p-1)}\triangleright \epsilon^i\hspace*{-3.5pt}&=&\hspace*{-3.5pt}
v^{ia_2\dots a_p}{\Bbb A}_{a_2\dots a_p}^{\;\;b_2\dots
b_p}e_{b_2}\otimes \dots 
\otimes e_{b_p}
\equiv v^{i\{a\}}{\Bbb A}_{\{a\}}^{\;\{b\}}\mathbf{e}_{\{b\}}\\
&=&\hspace*{-3.5pt}
q^pp_q\,v^{i\{a\}}\,\mathbf{e}_{\{b\}}\,v^{m\{b\}}\,B_m^{\, n}\,
u_{n\{a\}} = q^pp_q\,\epsilon^mB_m^{\,n}u_{n\{a\}}v^{i\{a\}}\\
&=&\hspace*{-3.5pt}q^{2p}\,\epsilon^m\,B_m^{\,n}C_n^{\,i} = \epsilon^i.
\end{eqnarray*}
Here at the last step we have used (\ref{nonsing}).

Therefore, under the action of the $q$-antisymmetrizer $A^{(p-1)}$ the
space 
$W={\rm Span}\{\epsilon^i\}$ is an invariant subspace in $V_{[p-1]}$ and
hence 
$W=V_{[p-1]}$. But as was proved in \cite{Gur} 
$$
\dim V_{[p-1]} = \dim V = n.
$$ 
Therefore the $n$ vectors $\epsilon^i$ cannot be linear dependent since
otherwise 
$\dim V_{[p-1]} = \dim W < n$. So, the set of eigenvector $\epsilon^i$
of the 
$q$-antisymmetrizer $A^{(p-1)}$ can be taken a basis of $V_{[p-1]}$.
\hfill\rule{6.5pt}{6.5pt}
\medskip

Now we are ready to establish the connection of $B$ and $R$ type
fundamental modules 
in the case of $R$-matrix with a finite symmetry rank.

\begin{utv}\label{u6}
Let the $R$-matrix has the symmetry rank $p$. Then the ${\cal SL}_q$
representation 
$\bar \theta_1$ obtained from (\ref{sl-r2}) is equivalent to
$\bar\pi_{[p-1]}$ obtained 
from (\ref{chi-a}) by sl-reduction (\ref{sl-nepr-q}).
\end{utv}

\noindent{\bf Proof\ \ }
Let us first consider the \Lq{} representation $\pi_{[p-1]}$
(\ref{chi-a}) which acts in 
the subspace $V_{[p-1]}$. In virtue of Lemma \ref{l3} we shall find the
matrices of 
operators $\pi_{[p-1]}(l_i^{\,j})$ in the basis of vectors $\epsilon^k$
(\ref{dual-bas}).
Using (\ref{chi-a}) we obtain (in the same notations as in the proof of
Lemma \ref{l3})
\begin{eqnarray*}
\frac{q^{2-p}}{(p-1)_q}\pi_{[p-1]}(l_i^{\,j})\triangleright\epsilon^k
\hspace*{-2.5mm}&=&\hspace*{-2.5mm}
v^{k\{a\}}{\Bbb A}_{\{a\}}^{\;m\{c\}}B_m^{\,j}{\Bbb
A}_{i\{c\}}^{\;\{b\}}
\mathbf{e}_{\{b\}} = ({\rm use\ } (\ref{*}))\\
&=&\hspace*{-2.5mm}q^{2p}p_q^2\,\mathbf{e}_{\{b\}}v^{r\{b\}}
B_s^{\,l}(v^{k\{a\}}u_{l\{a\}})B_r^{\,n}(v^{sm\{c\}}u_{ni\{c\}})B_m^{\,j}\\
&=&\hspace*{-2.5mm}q^{3p}p_q\,\epsilon^r\,(B_s^{\,l}C_l^{\,k})B_r^{\,n}
(v^{sm\{c\}}u_{ni\{c\}}) B_m^{\,j}=q^pp_q\,\epsilon^r
B_r^{\,n}(v^{km\{c\}}
u_{ni\{c\}})B_m^{\,j}.
\end{eqnarray*}
Introduce an $n^2\times n^2$ matrix $\Omega$ with matrix elements
$$
\Omega_{s_1s_2}^{\;r_1r_2}=
p_q(p-1)_q\sum_{\{a\}}u_{s_1s_2\{a\}}v^{r_1r_2\{a\}}.
$$
Then, the matrix of the operator $\pi_{[p-1]}(l_i^{\,j})$ in the basis
$\epsilon^k$ has the 
form (in compact notations)
\be
(\pi_{[p-1]}(L_2))_1 = q^{2(p-1)}\,B_1\Omega_{12}B_2.
\label{pi-matr}
\ee
With the use of (\ref{A-pr}) for $A^{(p)}$ and (\ref{C-uv}) for $C$ one
can express 
the matrix $\Omega$ in a more explicit form. Omitting straightforward
calculations we 
write down the final result
$$
\Omega_{12} = q^{2p-1}\Bigl(C_1C_2 - q C_1\Psi_{21}C_1\Bigr),
$$
where $\Psi$ is the skew-inverse to $R$-matrix as defined in
(\ref{closed}). Substituting
this into (\ref{pi-matr}) we find
$$
(\pi_{[p-1]}(L_2))_1 = q^{-3}I_{12} - q^{2p-2}\Psi_{21}C_1B_2.
$$
After sl-reduction (\ref{sl-nepr-q}) we get the representation of the
${\cal SL}_q$ algebra
$$
(\bar\pi_{[p-1]}(F_2))_1 = \frac{q^{1-p}}{(p-1)_q+q^{p+2}}\,
(I_{12} - q^{3p}p_q\Psi_{21}C_1B_2).
$$
The sl-reduction of the $R$ type representation (\ref{sl-r2}) leads in
turn to the result
$$
\bar\theta_1(F_2) = \frac{q^{p+1}}{(p-1)_q+q^{p+2}}\,(I_{12} -
q^{-p}p_q\,R_{12}).
$$
Next we take into account the connection of $\Psi$ and $R$ (see
Appendix)
\be
q^{2p}C_1\Psi_{21}B_2 = R^{-1}_{12}.
\label{psi-R}
\ee
With this formula one immediately gets
$$
C_1\,(\bar\pi_{[p-1]}(F_2))_1\,C^{-1}_1 = \bar\theta_1(F_2)
$$
which means that the corresponding modules are equivalent.\hfill
\rule{6.5pt}{6.5pt}

\subsection{Indecomposable modules: an example}

One can put a natural question about the completeness of the set of
representations thus 
obtained. In other words, whether an arbitrary finite dimensional REA
module with a non-commutative
representation of (\ref{RE}) is equivalent to a direct sum of modules
$V_\nu$ defined
in Proposition \ref{u2}? 

The answer to this question is negative. The matter is that REA
(\ref{RE}) possesses reducible
finite dimensional but {\it indecomposable} modules. The corresponding
mREA representations
do not admit a finite classical limit $q\rightarrow 1$. Let us give a
simplest example of such a 
module for algebra \Lq{} (\ref{ex:rea}). 

We start from the one-dimensional representation $\rho:\,\Lq\rightarrow
{\Bbb C}$
(see \cite{KulShweib})
\be
\rho(\hat L) = 
\left(
\matrix{
0&x\cr
y&z}
\right),\qquad x,y,z\in{\Bbb C}.
\label{one-d}
\ee
Then we use the comodule property (\ref{tsl}) in order to get the higher
dimensional
representation of \Lq. For this purpose take the known $R$-matrix
representation $\gamma$
of (\ref{RTT})
$$
\gamma(T_1) = P_{12}R_{12},\quad \gamma(S(T_1)) = R_{12}^{-1}P_{12}.
$$
Here $P$ is the transposition matrix, $R$ is the $U_q(sl_2)$ $R$-matrix
and the 
second matrix space stands for the representation space $V$, $\dim V =
2$. Then in accordance 
with (\ref{tsl}) we construct a two dimensional representation $\rho_2$
of \Lq{} in the space 
$V\otimes {\Bbb C} \cong V$
$$
\rho_2(\hat L_1)  = \gamma(T_1)\rho(\hat L_1)\gamma(S(T_1)) =
R_{21}\rho(\hat L_2)
R_{21}^{-1}.
$$
For the \Lq{} generators (\ref{ex:rea}) the explicit form of the
representation $\rho_2$
reads as follows
$$
\rho_2(\hat a) = \left(
\matrix{
0&-q\lambda x\cr
0&0}
\right),\;\;
\rho_2(\hat b) = \left(
\matrix{
qx &0\cr
0&q^{-1}x}
\right),\;\;
\rho_2(\hat c) = \left(
\matrix{
q^{-1}y&-\lambda z\cr
0&qy}
\right),\;\;
\rho_2(\hat d) = \left(
\matrix{
z&q^{-1}\lambda x\cr
0&z}
\right).
$$
The module $V$ with the representation $\rho_2$ is reducible. The
one-dimensional submodule
is spanned by the basis vector $e_1$. The corresponding one-dimensional
representation 
$$
\hat a\rightarrow 0, \quad \hat b\rightarrow qx,\quad
\hat c\rightarrow q^{-1}y,\quad \hat d\rightarrow z
$$
is connected with the initial one (\ref{one-d}) by an automorphism
$\eta$ of \Lq{} 
\cite{KulShweib}
$$
\eta\left(
\matrix{
\hat a&\hat b\cr
\hat c&\hat d}\right) = 
\left(
\matrix{
\hat a&\omega\hat b\cr
\omega^{-1}\hat c& \hat d
}\right),\qquad \forall\,\omega\in{\Bbb C}^{\times}.
$$
Nevertheless, being reducible, the module $V$ is obviously
indecomposable since matrices 
$\rho_2(\hat a)$ and $\rho_2(\hat d)$ cannot be transformed into
diagonal form (unless $x=0$). 
Therefore, this module cannot be presented as a direct sum of modules
$V_\nu$ constructed in 
Section \ref{th}. 

So, examining the completeness of the set of $V_\nu$ we have to reduce
the class 
of admissible modules to completely reducible ones and reformulate the
question
in the following way: is any {\it completely reducible} finite
dimensional module
over REA (\ref{RE}) isomorphic to a direct sum of modules $V_\nu$?

For an arbitrary $R$-matrix with a finite symmetry rank we have no
definite answer to 
this question. Given the only symmetry rank of $R$, one has too little
information
on the concrete structure of the corresponding REA. Perhaps, an analysis
of the explicit
commutation relations is needed here. The question on irreducibility of
modules $V_\nu$
themselves is also open in this case.

As for the $R$-matrix originated from the quantum universal enveloping
algebra
$U_q(sl_n)$ ($p=n$) it is highly plausible that the finite direct sums
of the modules
$V_\nu$ do exhaust all finite dimensional completely reducible
(non-commutative) 
representations of REA. The matter is that the matrix elements of the
corresponding 
representations are rational functions in $q$ with nonsingular limit
$q\rightarrow 1$. 
At that limit the mREA ${\cal SL}_q$ (\ref{sl-red}) tends to the algebra
$U(sl_n)$ and 
all the ${\cal SL}_q$ modules $V_\nu$ go to the corresponding modules
over $U(sl_n)$.
In particular the modules $V_\nu$ described in Proposition \ref{u2} must
be irreducible.

To conclude, we shortly summarize the main results and discuss some open 
problems and perspectives.

For the reflection equation algebra we have constructed the series of
finite
dimensional non-commutative representations which are parametrized by
Young diagrams.
The representations exist for any $R$-matrix satisfying the additional
conditions 
(\ref{Hec}), (\ref{closed}) and (\ref{f-rank}). The corresponding
modules $V_\nu$
are simple objects of a quasitensor Schur-Weyl category described in
detail in \cite{GLS1}.
As was pointed out in Section \ref{th}, the Grothendiek ring of the
Schur-Weyl category 
for the Hecke $R$-matrix with the symmetry rank $p$ is isomorphic to
that of the category 
of finite dimensional modules over $U(sl_p)$. Nevertheless, dimensions
of the modules and 
the characters of central elements could drastically differ from each
other.

Also, it is worth mentioning some further problems in this approach.
First of
them is the problem of constructing the representation theory for the
REA connected
with $R$-matrices of the Birman-Murakami-Wenzl type. Examples are given
by $R$-matrices
originated from the quantum groups of $B$, $C$ and $D$ series. The key
point here is
to develop the adequate technique for the $q$-analogs of the Young
idempotents.

Another interesting problem is the representation theory for the REA
with a
spectral parameter. This can find a lot of applications to the theory of
integrable systems.

%\newpage
\section*{Appendix}

This is a technical section where some auxiliary formulae of the main
text are proved.
First, we prove the trace formulae which were used in Proposition
\ref{u-sec}.
The decisive role belongs to the following result.

\begin{lem} Let $R$ be a solution of the Yang-Baxter equation
(\ref{YBE}), satisfying
the additional condition (\ref{closed}). Then
\be
Tr_{(0)}B_0R_{01}R_{02}^{-1} = P_{12}B_1,
\label{bas-ch}
\ee
where $P$ is the transposition matrix and B is defined in (\ref{BC}).
\end{lem}

\noindent{\bf Proof\ \ } Rewrite the Yang-Baxter equation (\ref{YBE}) in
the equivalent form
$$
R_{12}R_{23}R_{12}^{-1} = R_{23}^{-1}R_{12}R_{23}.
$$
Using this equation and definition (\ref{closed}) of the skew-inverse
matrix 
$\Psi$ we obtain the following relation
$$
Tr_{(0)}\Psi_{10}R_{02}R_{03}^{-1} =  P_{12}Tr_{(0)}R_{10}^{-1}R_{20}
\Psi_{03}P_{23}.
$$
Calculate now the trace in the first space. Since $Tr_{(1)}\Psi_{10} =
B_0$,
we get
$$
Tr_{(0)}B_0R_{02}R^{-1}_{03} =
Tr_{(01)}R_{20}^{-1}P_{12}R_{20}\Psi_{03}P_{23} =
Tr_{(0)}\Psi_{03}P_{23} = B_3P_{23} = P_{23}B_2.
$$
This result differs from (\ref{bas-ch}) only in the notations of the
matrix spaces.
\hfill\rule{6.5pt}{6.5pt}
\medskip

So, we are ready to prove the trace formulae used in Proposition
\ref{u-sec}.

{\bf i)}. ${\cal T}(n,k-1)\equiv Tr_{(1)}R^{-1}_{(1\rightarrow
n)}R^{-1}_{(k-1\rightarrow 1)}
B_1R_{1\,k+2}$ at $n<k-1$. First of all, one should use (\ref{R-ch}) in
order to draw the chain
$R^{-1}_{(1\rightarrow n)}$ to the right of $R^{-1}_{(k-1\rightarrow
1)}$
$$
R^{-1}_{(1\rightarrow n)}R^{-1}_{(k-1\rightarrow 1)} = 
R^{-1}_{(k-1\rightarrow 1)}R^{-1}_{(2\rightarrow n+1)}.
$$
The chain $R^{-1}_{(2\rightarrow n+1)}$ is evidently commute with
$B_1R_{1\,k+2}$,
therefore
$$
{\cal T}(n,k-1) = R^{-1}_{(k-1\rightarrow
2)}\left[Tr_{(1)}R^{-1}_{12}B_1R_{1\,k+2}
\right]R^{-1}_{(2\rightarrow n+1)}.
$$
Using the cyclic property of trace and then relation (\ref{bas-ch}), we
come to the desired 
result
$$
{\cal T}(n,k-1) = R^{-1}_{(k-1\rightarrow 2)}P_{2\,k+2}B_{k+2}
R^{-1}_{(2\rightarrow n+1)}.
$$
\medskip

{\bf ii)} The calculation of ${\cal T}(k-1,k-1)$ is more cumbersome. The
main difficulty
is that in this case the chains of $R$-matrices cannot be drawn through
each other.
As a consequence, it is not so easy to decrease the number of
$R$-matrices with the indices 
in the first space in order to apply (\ref{bas-ch}). However, with the
help of the Hecke 
condition and the Yang-Baxter equation the product of $R$-matrix chains
contained
in ${\cal T}(k-1,k-1)$ can be transformed as follows
$$
R^{-1}_{(1\rightarrow k-1)}R^{-1}_{(k-1\rightarrow 1)} = 
I_{12\dots k} - \lambda R_1^{-1} -
\lambda\sum_{n=2}^{k-1}R_{(n\rightarrow 2)}^{-1}
R_{1}^{-1}R_{(2\rightarrow n)}^{-1}.
$$
The terms in the right hand side contain at most one $R$ matrix with
indices in
the first space and hence, upon multiplying by $B_1R_{1\,k+2}$, we can
calculate
$Tr_{(1)}$ with the help of (\ref{bas-ch}). As a result we get the
formula which
was used in the proof of Proposition \ref{u-sec}
$$
{\cal T}(k-1,k-1) = I_{12\dots k} -\lambda
P_{2\,k+2}B_{k+2}-\lambda\sum_{n=2}^{k-1}
R^{-1}_{(n\rightarrow 2)}P_{2\,k+2}B_{k+2} R^{-1}_{(2\rightarrow n)}.
$$
\medskip

{\bf iii)}
At last, relation (\ref{psi-R}) is a direct consequence of
(\ref{bas-ch}). Indeed, 
multiply (\ref{bas-ch}) by $\Psi_{13}$ from the right and take the trace
in the first
space. Due to definition (\ref{closed}) of the matrix $\Psi$ we find
$$
Tr_{(0)}B_0P_{03}R_{02}^{-1} = Tr_{(1)}P_{12}B_1\Psi_{13} =
B_2\Psi_{23}\,Tr_{(1)}P_{12}
= B_2\Psi_{23}.
$$
On the other hand, due to the cyclic property of trace
$$
Tr_{(0)}B_0P_{03}R_{02}^{-1} = Tr_{(0)}P_{03}R_{02}^{-1}B_0 = 
R_{32}^{-1}B_3\,Tr_{(0)}P_{03} = R_{32}^{-1}B_3.
$$
Therefore
$$
B_2\Psi_{23} = R_{32}^{-1}B_3.
$$
Multiplying this by $C_3$ from the right and using (\ref{nonsing}) we
come
to the relation
$$
q^{2p}B_2\Psi_{23}C_3 = R_{32}^{-1.}
$$
Actually this is equivalent to (\ref{psi-R}), since basing on
(\ref{RBB}) one
can easily show that
$$
B_2\Psi_{23}C_3 = C_3\Psi_{23}B_2.
$$

\end{document}